\input amstex
\documentstyle {amsppt}
\magnification=1200
\nologo
\vsize=8.9truein
\hsize=6.5truein
\nopagenumbers

\def\abs#1{\left\vert #1 \right\vert}

\topmatter

\title
D\'ecompositions~de~groupes~par~produit~direct
\\
et groupes de Coxeter
\endtitle

\rightheadtext{D\'ecompositions~de~groupes~par~produit~direct}

\author
Yves de Cornulier et Pierre de la Harpe
\endauthor

\address
Yves de Cornulier, IGAT, EPFL, CH-1015 Lausanne.
Mel : decornul\@clipper.ens.fr
\endaddress

\address
Pierre de la Harpe, Section de Math\'ematiques, Universit\'e de Gen\`eve,
C.P. 64, 
\newline
CH-1211 Gen\`eve 4.
Mel : Pierre.delaHarpe\@math.unige.ch
\endaddress


\thanks
Nous remercions le Fonds National Suisse de la Recherche
Scientifique pour son soutien.
\endthanks

\abstract We provide examples of groups
which are indecomposable by direct product,
and more generally which are uniquely decomposable in direct products
of indecomposable groups.
Examples include Coxeter groups, for which we give an alternative approach
to recent results of L.~Paris.

For a finitely generated linear group $\Gamma$, we establish an upper bound
on the number of factors of which $\Gamma$ can be the direct product.
If $\Gamma$ is moreover torsion-free, it follows that $\Gamma$
is uniquely decomposable in direct product of indecomposable groups.

\medskip

\noindent
R\'ESUM\'E.
Nous montrons  des exemples de groupes
ind\'ecomposables par produit direct,
et plus g\'en\'eralement uniquement d\'ecomposables
en produits de groupes ind\'ecom\-posables.
Les exemples consid\'er\'es incluent les groupes de Coxeter,
pour lesquels nous red\'e\-montrons des r\'esultats r\'ecents de L.~Paris.

Pour un groupe lin\'eaire de type fini $\Gamma$, nous \'etablissons une borne
sup\'erieure sur le nombre de facteurs dont $\Gamma$ puisse \^etre produit direct.
Si $\Gamma$ est de plus sans torsion, il en r\'esulte que $\Gamma$
est uniquement d\'ecomposable
en produit de groupes ind\'ecom\-posables.

\smallskip

1. Introduction.

2. Premiers exemples.

3. D\'ecompositions de sous-groupes de groupes alg\'ebriques.

4. Groupes de Coxeter et groupes d'Artin.

5. Variation sur la notion d'ind\'ecomposabilit\'e.

6. Conditions Max et Min.

7. Groupes \`a quotients major\'es.

8. Majorations de nombres de facteurs directs.

\endabstract

\endtopmatter

\document

\head{\bf
1.~Introduction
}\endhead

Il existe plusieurs proc\'ed\'es pour \lq\lq d\'ecomposer\rq\rq \ des groupes
en constituants qu'on peut imaginer d'\'etude \lq\lq plus simple\rq\rq .
L'objet du pr\'esent article est la d\'ecomposition par {\it produit direct.}

\medskip

Un groupe est {\it ind\'ecomposable} 
s'il n'est pas un produit direct de mani\`ere non banale.
Notre premier but est d'\'etablir l'ind\'ecomposabilit\'e
de certains groupes apparaissant dans des contextes g\'eom\'etriques~; 
par exemple des sous-groupes \lq\lq assez grands\rq\rq \ 
ou m\^eme des r\'eseaux 
dans certains groupes de Lie ou groupes alg\'ebriques.

\medskip

La condition pour un groupe
d'\^etre isomorphe \`a un produit fini de groupes ind\'ecompo\-sables 
n'est pas tr\`es restrictive,
et il est bien connu que de nombreux groupes
poss\`edent plusieurs telles d\'ecompositions 
qui sont essentiellement diff\'erentes.
Voir l'exemple $A_{2,2} \times A_{3,3}$ de \cite{Kuro, \S~42}
et les exemples ab\'eliens rappel\'es ci-dessous
(num\'eros (xiv) et (xvii) du chapitre 2)~;
il y a aussi des exemples \lq\lq plus dramatiques\rq\rq \ 
dus \`a  Baumslag~:
pour toute paire d'entiers $m,n \ge 2$,
il existe des groupes de type fini $A_1, \hdots, A_m, B_1, \hdots, B_n$, 
nilpotents et sans torsion, ind\'ecomposables et non isomorphes deux \`a deux,
tels que les produits directs $A_1 \times \cdots \times A_m$
et $B_1 \times \cdots \times B_n$ sont isomorphes \cite{Baum}.
Il existe de nombreux autres cas de \lq\lq mauvais comportements\rq\rq \
du point de vue des produits directs~: 
par exemple des groupes $\Gamma$ de type fini 
\footnote{
C'est un probl\`eme ouvert bien \'etabli de savoir s'il existe
un groupe $\Gamma$ infini de pr\'esentation finie 
qui soit isomorphe \`a $\Gamma \times \Gamma$.
Voici une question peut-\^etre moins ambitieuse~: 
existe-t-il un groupe $\Gamma$ infini de pr\'esentation finie
tel que, pour tout entier $m_0 \ge 1$,
il existe des groupes $\Gamma_1, \hdots, \Gamma_m$ ($m \ge m_0$)
non r\'eduits \`a un \'el\'ement 
dont le produit direct soit isomorphe \`a $\Gamma$~? 
}
isomorphes \`a
$\Gamma \times \Delta$, $\Delta \ne \{1\}$ \cite{Hi86}
ou m\^eme \`a $\Gamma \times \Gamma$ (\cite{Ty74}, \cite{Ty80}, \cite{Meie}).

La situation est donc beaucoup moins simple 
que celle qui pr\'evaut pour les produits libres,
puisque le th\'eor\`eme de Grushko assure que  tout groupe de type fini 
est  produit libre d'un nombre fini de groupes librement ind\'ecomposables, 
et ceci de mani\`ere essentiellement unique
(voir par exemple \cite{Stal} ou \cite{ScWa}).

La plupart des complications qui apparaissent dans les d\'ecompositions 
d'un groupe~$G$ par produit direct sont li\'ees au centre de $G$. 
En effet, si celui-ci est r\'eduit \`a un \'el\'ement, 
une d\'ecomposition de $G$ en produit fini de groupes ind\'ecomposables
est n\'ecessairement unique~; ce fait, bien connu, est rappel\'e
ci-dessous comme cons\'equence de la proposition~2.

\medskip

Convenons qu'un groupe est {\it uniquement directement d\'ecomposable}
\footnote{
D'autres auteurs \cite{Hi90} \'ecrivent \lq\lq groupe R.K.S.\rq\rq ,
en r\'ef\'erence \`a Remak, Krull et Schmidt.
Pour une notion de d\'ecomposabilit\'e apparent\'ee mais distincte,
voir ci-dessous le chapitre 5.
}
s'il est somme res\-treinte de groupes ind\'ecomposables,
et ceci de mani\`ere unique \`a isomorphisme pr\`es des facteurs.
(Rappelons que la {\it somme restreinte} 
d'une famille $(\Gamma_{\iota})_{\iota \in I}$  d\'esigne 
le sous-groupe du groupe produit des $G_{\iota}$ form\'e des \'el\'ements
$(\gamma_{\iota})_{\iota \in I}$ tels que 
$\gamma_{\iota} = 1$ pour presque tout $\iota$~;
dans le cas o\`u $I$ est fini, 
nous suivons Bourbaki en \'ecrivant aussi \lq\lq produit direct\rq\rq \ 
au lieu de \lq\lq somme restreinte\rq\rq .)
Remarquons qu'un groupe uniquement directement d\'ecomposable 
qui est de type fini est n\'ecessairement
produit direct d'un nombre {\it fini} de facteurs ind\'ecomposables.

Notre second but est d'\'etablir que certaines classes de groupes
sont uniquement directement d\'ecomposables (propositions 3, 5 et 14).
Apr\`es des rappels d'exemples standard, 
nous analysons en particulier la d\'ecomposabilit\'e des
{\it sous-groupes Zariski-denses de groupes alg\'ebriques}.

Nous illustrons notre m\'ethode par les
{\it groupes de Coxeter} (proposition 8).
Dans \cite{Par2}, Paris a montr\'e 
de mani\`ere compl\`etement diff\'erente
que les groupes de Coxeter de type fini 
sont uniquement directement d\'ecomposables.
L'origine du pr\'esent travail est  
la recherche d'une autre preuve de ce fait.

\medskip

Enfin, nous \'etablissons une borne sup\'erieure 
pour le nombre de facteurs non r\'eduits \`a un \'el\'ement 
dans une d\'ecomposition par produit direct 
d'un groupe lin\'eaire de type fini
(th\'eor\`eme 13).

\medskip

La d\'ecomposition des groupes par produit direct intervient dans des
r\'esultats de d\'ecom\-position d'espaces topologiques. 
Par exemple, soit $\Gamma$ le groupe fondamental d'une vari\'et\'e riemannienne 
compacte \`a courbure n\'egative ou nulle~;
supposons le centre de $\Gamma$ r\'eduit \`a un \'el\'ement.
C'est un cas particulier de
r\'esultats classiques de Gromoll-Wolf (1971) et Lawson-Yau (1972)
qu'une d\'ecomposition du groupe 
en produit direct $\Gamma = \Gamma_1 \times \Gamma_2$
correspond n\'ecessai\-rement \`a une d\'ecomposition de la vari\'et\'e
en produit riemannien.
Voici un corollaire de r\'esultats plus r\'ecents
(voir \cite{BrHa}, chapitre II.6, ainsi que \cite{Schr}
et le chapitre 10 de \cite{Eber}).

\medskip

{\it Soit $Y$ un espace g\'eod\'esique compact 
qui poss\`ede  la propri\'et\'e d'extension des g\'eod\'e\-siques
et qui est \`a courbure n\'egative ou nulle.
Supposons que le groupe fondamental de $Y$ est un produit direct
$\Gamma = \Gamma_1 \times \Gamma_2$ 
et que le centre de $\Gamma$ est r\'eduit \`a un \'el\'ement.
Alors $Y$ est un produit d'espaces m\'etriques $Y_1 \times Y_2$ 
de telle sorte que le groupe fondamental de $Y_i$ est $\Gamma_i$ ($i=1,2$).
}

\medskip

Mentionnons encore que, en th\'eorie ergodique, 
certains produits directs donnent lieu 
\`a des ph\'enom\`enes de rigidit\'e remarquables \cite{MoSh}.

\bigskip
\head{\bf
2.~Premiers exemples
}\endhead

Un groupe $\Gamma$  est {\it ind\'ecomposable} si 
\footnote{
Nous adoptons ici une terminologie selon laquelle
le groupe $\{1\}$ n'est ni ind\'ecomposable, ni d\'ecompo\-sable
(mais n\'eanmoins le produit de la famille vide de groupes ind\'ecomposables).
Ceci est coh\'erent avec la convention (par exemple de Bourbaki)
selon laquelle le groupe \`a un \'el\'ement n'est pas simple~!
}
$\Gamma \ne \{1\}$ et si,
pour tout isomorphisme de $\Gamma$ 
avec un produit direct $\Gamma_1 \times \Gamma_2$,
l'un des groupes $\Gamma_1,\Gamma_2$ est r\'eduit \`a un \'el\'ement.
La notion est classique~: voir parmi d'autres 
\cite{Kuro, \S \ 17},
\cite{Rotm, chap. 6} et \cite{Suzu, \S \ 2.4}.

\bigskip

   (i) Tout groupe simple est ind\'ecomposable. 

\medskip

   (ii) Il existe de nombreux groupes finis ind\'ecomposables 
qui ne sont pas simples.
Par exemple, pour tout $n \ge 2$,
le groupe sym\'etrique $\operatorname{Sym}_n$, 
qui est un groupe de Coxeter de type $A_n$,
est ind\'ecomposable
(cela r\'esulte
de ce que tout sous-groupe normal distinct de~$\{1\}$ 
contient le groupe altern\'e simple $\operatorname{Alt}_n$ si $n \ne 2,4$,
et d'un argument direct si $n=2$ ou $4$)~; 
pour tout $n \ge 3$, 
le produit semi-direct standard 
$\operatorname{Sym}_n \ltimes (\bold Z / 2\bold Z)^{n-1}$,
qui est un groupe de Coxeter de type $D_n$, 
est ind\'ecomposable (l'argument de la proposition 1 ci-dessous s'applique).
En revanche, on sait qu'il existe des groupes de Coxeter irr\'eductibles finis 
sont d\'ecomposables,
et alors d\'ecomposables en produit direct d'un groupe ind\'ecomposable 
et du centre \`a deux \'el\'ements~;
ce sont les groupes di\'edraux d'ordres $8n+4$
et les groupes de Coxeter de type type $B_n$ pour $n \ge 3$ impair
(la v\'erification  est laiss\'ee au lecteur),
ainsi que les groupes de type $H_3$ et $E_7$
(voir \cite{BouL}, chapitre 6, \S~4, exercices 3 et 11).

Pour tout entier $n$ pair,   
il existe un groupe fini m\'etab\'elien ind\'ecomposable d'ordre~$n$~;
par exemple le produit semi-direct 
$$
(\bold Z / 2^k \bold Z) \ltimes_{\pm} (\bold Z / m\bold Z) 
\, = \,
\langle a,b 
\hskip.1cm \vert \hskip.1cm 
a^{2^k}=1, b^m = 1, a^{-1}ba = b^{-1} \rangle ,
$$
o\`u $k,m \ge 1$ sont tels que $n = 2^km$ avec $m$ impair.
En revanche, les entiers impairs 
qui sont des ordres de groupes finis ind\'ecomposables 
constituent un ensemble d'entiers de densit\'e nulle \cite{ErPa}.

\medskip

   (iii) Les groupes ab\'eliens $\bold Z$ et $\bold Q$ sont ind\'ecomposables~;
plus g\'en\'eralement, 
tout sous-groupe de $\bold Q$ non r\'eduit \`a un \'el\'ement
est ind\'ecomposable (pour la description de ces groupes,
voir par exemple le chapitre 10 de \cite{Rotm}).
Soient $p$ un nombre premier et $m \ge 1$ un entier~;
un groupe cyclique d'ordre $p^m$ est ind\'ecomposable,
de m\^eme que le sous-groupe $\bold Z(p^{\infty})$ de $\bold C^*$ 
des racines de l'unit\'e d'ordres des puissances de $p$.
En effet, dans chacun de ces groupes, 
l'intersection de deux sous-groupes non r\'eduits \`a un \'el\'ement
n'est jamais r\'eduite \`a un \'el\'ement.
Le groupe additif $\bold Z_p$ des entiers $p$-adiques 
(vu comme groupe discret) est ind\'ecomposable
\cite{Kapl, Section 15}.

Un groupe ab\'elien divisible ind\'ecomposable est isomorphe \`a l'un de
$\bold Q$, $\bold Z(p^{\infty})$~; les groupes ab\'eliens divisibles
sont uniquement directement d\'ecomposables \cite{Fuch, \S~23}.
Un groupe ab\'elien ind\'ecomposable est ou bien de torsion ou bien sans torsion~;
s'il est de torsion, il est isomorphe \`a l'un de
$\bold Z /p^m\bold Z$, $\bold Z(p^{\infty})$~;
voir \cite{Kapl, Section 9}.
Les groupes ab\'eliens sans torsion sont beaucoup moins bien compris,
et certainement pas tous uniquement directement d\'ecomposables~;
voir (xvii) ci-dessous.

\medskip

   (iv) Un groupe nilpotent $\Gamma$ dont le centre est ind\'ecomposable
est lui-m\^eme ind\'ecom\-posable.
En effet, soit $\Gamma = \Gamma_1 \times \Gamma_2$ 
une d\'ecomposition en produit direct.
Le centre de $\Gamma$ \'etant \'egal 
au produit des centres de $\Gamma_1$ et $\Gamma_2$,
l'un de ceux-ci est r\'eduit \`a un \'el\'ement.
L'assertion r\'esulte de ce qu'un
groupe nilpotent dont le centre est r\'eduit \`a un \'el\'ement
est lui-m\^eme  r\'eduit \`a un \'el\'ement.

En particulier, le {\it groupe de Heisenberg}
$\left( \matrix 
1 & \bold Z & \bold Z \\ 0 & 1 & \bold Z \\ 0 & 0 & 1 
\endmatrix \right)$
est ind\'ecomposable.

Notons 
$\Gamma = C^1(\Gamma) \supset \cdots \supset 
C^{j+1}(\Gamma) = [\Gamma,C^j(\Gamma)] \supset \cdots$
la s\'erie centrale descendante d'un groupe $\Gamma$.
Soient $k,j \ge 2$ des entiers et
$F_k$ le groupe non ab\'elien libre \`a $k$ g\'en\'erateurs~;
le {\it groupe nilpotent libre} de classe $j$ \`a $k$ g\'en\'erateurs
$\Gamma = F_k/C^{j+1}(F_k)$ est ind\'ecomposable.

En effet, soit $\Gamma = \Gamma_1 \times \Gamma_2$
une d\'ecomposition en produit direct. 
Soient $l,m$ les rangs des groupes ab\'eliens libres
$\Gamma_1/C^2(\Gamma_1), \Gamma_2/C^2(\Gamma_2)$, respectivement~;
notons que le rang~$k$ du groupe ab\'elien libre 
$\Gamma/ C^2(\Gamma) = F_k/C^2(F_k)$
est \'egal \`a la somme $l+m$.
D'une part, le rang du groupe ab\'elien libre 
$C^2(\Gamma) / C^3(\Gamma) = C^2(F_k) / C^3(F_k)$ 
est le coefficient binomial $\binom{k}{2}$.
D'autre part, le rang du groupe ab\'elien libre $C^2(\Gamma_1) / C^3(\Gamma_1)$
est major\'e par le rang $\binom{l}{2}$ de $C^2(F_l) / C^3(F_l)$,
et de m\^eme pour $\Gamma_2$ et $\binom{m}{2}$. Comme
$$
C^2(\Gamma) / C^3(\Gamma) \, \approx \,
\left(C^2(\Gamma_1) / C^3(\Gamma_1)\right) \times 
\left(C^2(\Gamma_2) / C^3(\Gamma_2)\right) ,
$$
il en r\'esulte que
$\binom{k}{2} \le \binom{l}{2} + \binom{m}{2}$,
donc que l'un de $k,l$ est z\'ero,
et par suite que l'un de $\Gamma_1,\Gamma_2$ est r\'eduit \`a un \'el\'ement.

\medskip

   (v) Un groupe r\'esiduellement r\'esoluble $\Gamma$ dont l'ab\'elianis\'e 
$\Gamma / [\Gamma,\Gamma]$
est ind\'ecomposable est lui-m\^eme ind\'ecomposable.
Par exemple, le groupe
$\left( \matrix 2^{\bold Z} & \bold Z [1/2] \\ 0 & 1 \endmatrix \right)$
est ind\'ecomposable.
(L'argument de la proposition 3 
permet de remplacer $2$ par un autre nombre rationnel.)

   Plus g\'en\'eralement
un groupe r\'esiduellement r\'esoluble $\Gamma \ne \{1\}$ 
est ind\'ecomposable d\`es qu'un 
des quotients $\Gamma / C^j(\Gamma)$ 
ou  $\Gamma / D^j(\Gamma)$ 
est ind\'ecomposable pour un entier $j\ge 2$, 
o\`u $\Gamma = D^0(\Gamma) \supset \cdots \supset 
D^{i+1}(\Gamma) = [D^i(\Gamma),D^i(\Gamma)] \supset \cdots$
d\'esigne la s\'erie d\'eriv\'ee.

D\'etaillons par exemple l'argument pour $C^j$~:
soit $\Gamma = \Gamma_1 \times \Gamma_2$ un groupe r\'esiduellement r\'esoluble
tel que $\Gamma/C^j(\Gamma)$ est ind\'ecomposable~;
nous pouvons supposer les notations telles que 
$\Gamma_2/C^j(\Gamma_2)$ est r\'eduit \`a un \'el\'ement~; 
{\it a fortiori,} $\Gamma_2/[\Gamma_2,\Gamma_2]$ est r\'eduit \`a un \'el\'ement,
ce qui implique 
que $\Gamma_2 = \{1\}$ puisque $\Gamma_2$ est r\'esiduellement r\'esoluble.

\medskip

   (vi) Le groupe fondamental
$\Gamma = \langle x,y \hskip.1cm \mid \hskip.1cm xyx^{-1} = y^{-1} \rangle$
d'une bouteille de Klein est ind\'ecomposable. 

En effet, soit $\Gamma = \Gamma_1 \times \Gamma_2$ 
une d\'ecomposition en produit direct.
Le groupe d\'eriv\'e $D\Gamma = D\Gamma_1 \times D\Gamma_2$ 
est engendr\'e par $y^2$~; il est cyclique infini, et donc ind\'ecomposable,
de sorte qu'on peut supposer $D\Gamma_2 = \{1\}$, 
c'est-\`a-dire $\Gamma_2$ ab\'elien.
Pour $j=1,2$, notons $x_j$ [respectivement $y_j$] la projection de $x$ [resp. $y$] 
sur $\Gamma_j$.
Alors $D\Gamma_2 = \{1\}$ implique $y_2^2 = 1$, 
et en fait $y_2 = 1$ puisque $\Gamma$ est sans torsion.
Le centre $Z(\Gamma) = Z(\Gamma_1) \times Z(\Gamma_2)$
est engendr\'e par $x^2$~; 
il est \'egalement cyclique infini, donc ind\'ecomposable.
Si nous avions $Z(\Gamma_1) = \{1\}$, 
nous aurions aussi $x_1 = 1$ comme ci-dessus,
donc $\Gamma = \langle y_1 \rangle \times \langle x_2 \rangle$ serait ab\'elien,
ce qui est absurde. 
C'est donc $Z(\Gamma_2)$ qui est r\'eduit \`a un \'el\'ement,
ce qui montre enfin que $\Gamma_2$ lui-m\^eme est r\'eduit \`a un \'el\'ement.

\medskip

   (vii) Un produit libre de deux groupes non r\'eduits \`a un \'el\'ement est
ind\'ecomposable.  
\footnote{
Notons la cons\'equence suivante pour le probl\`eme de d\'ecision
relatif \`a la d\'ecomposition en produit direct~:
le probl\`eme de savoir si un groupe donn\'e par une pr\'esentation
est d\'ecomposable ou non n'est pas algorithmiquement r\'esoluble. 
En effet, soient $A$ et $B$ deux groupes de pr\'esentation finie
non r\'eduits \`a un \'el\'ement
(par exemple $A = B = \bold Z$)
et $\Delta$ un groupe de pr\'esentation finie.
Le groupe $\Gamma = (A \times B)\ast \Delta$
est d\'ecomposable si et seulement si $\Delta$ est r\'eduit \`a un \'el\'ement,
et il est bien connu qu'il n'existe pas d'algorithme permettant
de savoir si un groupe donn\'e par une pr\'esentation finie 
est ou n'est pas r\'eduit \`a un \'el\'ement.
}
Plus g\'en\'eralement, un sous-groupe d'un produit libre $\Gamma_1 \ast \Gamma_2$
qui n'est conjugu\'e ni \`a un sous-groupe de $\Gamma_1$ ni \`a un sous-groupe de
$\Gamma_2$ est ind\'ecomposable. Voir \cite{Kuro}, \S \S~34 et 35.

\medskip

   (viii) Soit $\Gamma$ un sous-groupe non r\'eduit \`a un \'el\'ement
d'un groupe  qui est hyperbolique au sens de Gromov et sans torsion.
Alors $\Gamma$ est ind\'ecomposable car, pour tout \'el\'ement $\gamma \ne 1$
dans un tel groupe, le centralisateur de $\gamma$ est cyclique infini.

\medskip

   (ix) Soit $G$ un groupe de Lie r\'eel connexe semi-simple,
\`a centre fini et sans facteur compact, de rang r\'eel au moins deux,
et soit $\Gamma$ un r\'eseau irr\'eductible dans $G$.
On~suppose le centre de $\Gamma$ r\'eduit \`a un \'el\'ement.
C'est une cons\'equence imm\'ediate du th\'eor\`eme de Margulis concernant
les sous-groupes normaux de $\Gamma$ (voir le chapitre 8 de \cite{Zimm}) 
que $\Gamma$ est ind\'ecomposable~;
ceci s'applique par exemple \`a 
$\Gamma = \operatorname{PSL}_d(\bold Z)$ pour $d \ge 3$.
Avec des formulations plus g\'en\'erales (voir le chapitre VIII de \cite{Marg}),
on montre de m\^eme que des groupes comme $\operatorname{PSL}_d(\bold Z[1/p])$,
qui est un r\'eseau irr\'eductible dans 
$\operatorname{PSL}_d(\bold R) \times \operatorname{PSL}_d(\bold Q_p)$,
sont ind\'ecomposables ($p$ est un nombre premier, $d$ un entier, $d \ge 2$).
 
Le m\^eme type de r\'esultat (et d'argument) vaut encore plus g\'en\'eralement 
pour des r\'eseaux irr\'eductibles 
dans certains produits de groupes localement compacts \cite{BaSh}.

Dans de nombreux cas, ces affirmations peuvent \^etre d\'emontr\'ees
par des arguments plus \'economiques.
Par exemple, les groupes $\operatorname{PSL}_d(\bold Z[1/p])$ 
sont denses dans le groupe de Lie $\operatorname{PSL}_d(\bold R)$,
a fortiori Zariski-denses dans le groupe alg\'ebrique 
$\operatorname{PSL}_d(\bold C)$,  et sont donc ind\'ecomposables en vertu de la
proposition 3 ci-dessous.

\medskip

   (x) Soit $\Gamma$ un groupe tel que, 
pour toute paire $N_1,N_2$ de sous-groupes normaux
non r\'eduits \`a un \'el\'ement, 
l'intersection $N_1 \cap N_2$ ne l'est pas non plus~; 
alors $\Gamma$ est \'evidemment ind\'ecomposable.
Dans ce cas, les sous-groupes normaux non r\'eduits 
\`a un \'el\'ement forment une base de voisinage de $1$ 
pour une  topologie sur $\Gamma$
qu'on appelle la {\it topologie pro-normale} 
et qui est \'etudi\'ee dans \cite{GeGl}.

   Un sous-groupe Zariski-dense $\Gamma$ d'un groupe de Lie $G$ connexe simple
\`a centre trivial poss\`ede cette propri\'et\'e.
En effet, soient $N_1$ et $N_2$ deux sous-groupes non r\'eduits \`a un \'el\'ement
de $\Gamma$. 
D\'esignons par $\overline{N}_1$ et $\overline{N}_2$ leurs adh\'erences de Zariski~;
ce sont des sous-groupes normaux du groupe simple $G$, 
de sorte qu'ils co\"{\i}ncident avec $G$.
Si nous avions $N_1 \cap N_2 = \{1\}$, nous aurions aussi $[N_1,N_2] = \{1\}$
et $[\overline{N}_1, \overline{N}_2] = [G,G] = \{1\}$,
ce qui est absurde.

   Le \lq\lq groupe de Grigorchuk\rq\rq \ poss\`ede aussi cette propri\'et\'e.
C'est le $2$-groupe infini de type fini qui appara\^{\i}t dans \cite{Grig}~;
voir aussi, par exemple, le th\'eor\`eme VIII.42 de \cite{Harp}.
Ceci s'\'etend \`a tout groupe {\it juste infini}, c'est-\`a-dire \`a tout
groupe infini dont tous les quotients propres sont finis~;
voir \cite{SaSS}, en particulier le chapitre r\'edig\'e par J.~Wilson.

   De m\^eme pour le \lq\lq groupe $F$ de Thompson\rq\rq  ,
dont on sait que le groupe d\'eriv\'e est d'une part simple
et d'autre part contenu dans tout sous-groupe normal non r\'eduit \`a un \'el\'ement
(th\'eor\`eme 4.5 et preuve du th\'eor\`eme 4.3 dans \cite{CaFP}).

\medskip

   (xi) Nous d\'emontrons ci-dessous l'ind\'ecomposabilit\'e 
des groupes de Coxeter irr\'educ\-tibles infinis 
(corollaire  de la proposition 1 et proposition 8).

\medskip

   (xii)  Soit $\Gamma$ un groupe non r\'eduit \`a un \'el\'ement qui est
de pr\'esentation finie et de dimension cohomologique au plus 2.
Alors $\Gamma$  est ou bien ind\'ecomposable, ou bien un produit direct de deux
groupes libres.
En effet, s'il  existe deux sous-groupes $\Gamma _1,\Gamma_2$ de $\Gamma$ 
tels que $\Gamma = \Gamma_1 \times \Gamma_2$,
un r\'esultat de Bieri implique que 
$\Gamma_1$ et $\Gamma_2$ sont de dimension cohomologique au plus un
[Bier, corollaire 8.6],
de sorte que $\Gamma_1$ et $\Gamma_2$ sont libres
par un th\'eor\`eme de Stallings
(voir par exemple [Bier, theoreme 7.6]).

   Voici deux familles d'exemples de groupes
de dimension cohomologique au plus~2~:
les sous-groupes sans torsion des groupes \`a un relateur
([Bier, th\'eor\`eme 7.7], r\'esultat d\^u a Lyndon pour un
groupe \`a un relateur sans torsion),
et 
les groupes fondamentaux
des vari\'et\'es non compactes de dimension $3$,
en particulier les groupes de noeuds
(voir par exemple \cite{Serr}, no 1.5 et lemme 5 du no 2.1).

   Bagherzadeh a montr\'e un r\'esultat plus pr\'ecis~:
les seuls sous-groupes d\'ecomposables 
des groupes \`a un relateur sans torsion 
sont des produits directs d'un groupe cyclique infini et d'un groupe libre 
(voir le corollaire 4.9 de \cite{Bagh}).
En particulier, tout sous-groupe d\'ecomposable de type fini d'un
groupe \`a un relateur sans torsion est isomorphe \`a $\bold Z^2$.

\medskip

   (xiii) Pour toute paire $k,l$ d'entiers, $k,l \ge 2$, le produit amalgam\'e
$$
A_{k,l} \, = \, \langle a_1,a_2 \hskip.2cm \vert \hskip.2cm a_1^k = a_2^l \rangle
$$
est ind\'ecomposable~; 
r\'ep\'etons l'argument simple du livre de Kurosh.

   Le groupe $A_{k,l}$ est sans torsion, 
et son centre est cyclique infini engendr\'e par $a_1^k = a_2^l$
(ce sont des cons\'equences imm\'ediates de r\'esultats concernant les formes 
normales dans les produits amalgam\'es, voir par exemple le \S~4.2 de \cite{MaKS}). 
Si $A_{k,l} = X \times Y$ est une d\'ecomposition en produit direct
et si $(x,y) \in X \times Y$ 
est l'\'ecriture du g\'en\'erateur $a_1^k$ du centre de $A_{k,l}$,
il en r\'esulte que l'un de $x,y$ est $1$~;
sans restreindre la g\'en\'eralit\'e de ce qui suit, 
nous pouvons supposer que $y = 1$.
Notons $a_1 = (x_1,y_1)$, $a_2 = (x_2,y_2)$ 
les \'ecritures dans le produit $X \times Y$
des g\'en\'erateurs $a_1,a_2$ de $A_{k,l}$~;
remarquons que le groupe $X$ [respectivement le groupe $Y$]
est engendr\'e par $x_1$ et $x_2$ [resp. par $y_1$ et $y_2$].
Comme $Y$ est sans torsion et comme $y_1^k = y_2^l = y = 1$, 
nous avons $y_1 = y_2 = 1$.
Il en r\'esulte que $Y$ est r\'eduit \`a un \'el\'ement.

\medskip

   (xiv) Avec les notations de (xiii), 
consid\'erons un entier $k \ge 2$ et le produit direct
$$
\Gamma \, = \, 
A_{k,k} \times A_{k+1,k+1} \, = \, 
\langle a_1,a_2 \hskip.2cm \vert \hskip.2cm a_1^k = a_2^k \rangle
\, \times \,
\langle b_1,b_2 \hskip.2cm \vert \hskip.2cm b_1^{k+1} = b_2^{k+1} \rangle .
$$
Posons $a = a_1^k = a_2^k$, qui est un g\'en\'erateur du centre de $A_{k,k}$,
et $b = b_1^{k+1} = b_2^{k+1}$, qui est un g\'en\'erateur du centre de $A_{k+1,k+1}$.
Dans $\Gamma$, posons
$$
\aligned
c_1 \, &= \, ab^{-1}a_1 ,\quad
c_2 \, = \, ab^{-1}a_2 ,\quad
c_3 \, = \, ab^{-1}b_1 ,\quad
c_4 \, = \, ab^{-1}b_2  \\
c \, &= \, c_1^k \, = \, c_2^k \, = \, c_3^{k+1} \, = \, c_4^{k+1}, \qquad
d \, = \, ab^{-1} .
\endaligned
$$
On v\'erifie que le sous-groupe $C$ de $\Gamma$
engendr\'e par $c_1,c_2,c_3,c_4$ 
a un centre cyclique infini engendr\'e par $c$, 
et on montre comme en (xiii) que $C$ est ind\'ecomposable.

   L'int\'er\^et de cet exemple est le suivant~: 
le groupe $A_{k,k} \times A_{k+1,k+1}$
est isomorphe au produit direct du groupe $C$ et du groupe cyclique infini
$D$ engendr\'e par $d$~;
et les groupes $A_{k,k}, A_{k+1,k+1}, C, D$, tous ind\'ecomposables, 
sont non isomorphes deux \`a deux.
En particulier, le groupe de pr\'esentation finie
$A_{k,k} \times A_{k+1,k+1}$ n'est pas uniquement directement
d\'ecomposable. Tout ceci appara\^{\i}t dans \cite{Kuro, \S~42}, pour $k=2$.

   Cet exemple montre aussi qu'un r\'eseau 
dans un groupe de Lie r\'eel semi-simple
peut ne pas \^etre uniquement directement d\'ecomposable. 
En effet, pour $k \ge 3$, le quotient 
$\langle a_1,a_2 \hskip.2cm \vert \hskip.2cm a_1^k = a_2^k = 1 \rangle$
de $A_{k,k}$ par son centre 
est un r\'eseau dans le groupe de Lie $\operatorname{PSL}_2(\bold R)$~;
c'est un groupe fuchsien qui poss\`ede 
un quadrilat\`ere fondamental dans le demi-plan de Poincar\'e
ayant deux sommets oppos\'es d'angle $\pi/k$ 
et les deux autres  sommets \`a l'infini.
Le groupe $A_{k,k}$ lui-m\^eme est un r\'eseau dans
le rev\^etement universel $\tilde G$ de $\operatorname{PSL}_2(\bold R)$.
Par suite, $A_{k,k} \times A_{k+1,k+1}$
est un r\'eseau r\'eductible dans le groupe de Lie semi-simple
qui est produit direct de deux copies
de $\tilde G$.

Notons enfin que les groupes $A_{k,k} \times A_{k+1,k+1}$ sont lin\'eaires.
En effet, bien que le groupe $\tilde G$ ne soit pas lin\'eaire,
ses r\'eseaux $A_{k,k}$ le sont
(comparer avec la proposition de Toledo, Millson et Gersten
qui appara\^{\i}t au no IV.48 de \cite{Harp}).

\medskip

   Citons plus bri\`evement quelques exemples illustrant la notion d'unique
d\'ecomposabilit\'e directe.

\medskip

   (xv) Un groupe ind\'ecomposable 
est \'evidemment uniquement directement d\'ecompo\-sable. 
Avec nos conventions, le groupe \`a un \'el\'ement est uniquement
directement d\'ecom\-posable.

\medskip

   (xvi) Tout groupe fini est uniquement directement d\'ecomposable.
C'est une cons\'equence imm\'ediate 
du th\'eor\`eme de Wedderburn-Remak-Krull-Schmidt,
mais un r\'esul\-tat non banal~!
(Voir aussi ci-dessous la  proposition 9.)

\medskip

   (xvii)  Un groupe ab\'elien de type fini est uniquement directement
d\'ecomposable (voir par exemple \cite{BouA'}, chap. VII, 
\S \ 4, no 8).
Un groupe ab\'elien libre est uniquement directement d\'ecomposable~;
par exemple, le groupe multiplicatif $\bold Q_+^*$, 
isomorphe \`a la somme restreinte de copies de $\bold Z$ 
index\'ees par les nombres premiers, 
est uniquement directement d\'ecomposable. 
Une somme restreinte de groupes ab\'eliens de rang $1$, 
c'est-\`a-dire de sous-groupes de $\bold Q$, 
est uniquement directement d\'ecomposable.
(C'est un r\'esultat de Baer~; voir par exemple \cite{Fuch, Section 86}.)

\smallskip

   En revanche, les groupes ab\'eliens de rang fini sans torsion 
sont loins d'\^etre tous uniquement directement d\'ecomposables. 
Par exemple, pour tout entier $n \ge 2$,
il existe des groupes ab\'eliens de rang fini sans torsion
$A, A', A_1, \hdots, A_n$, 
ind\'ecomposables et non isomorphes deux \`a deux,
tels que $A \oplus A'$ et $A_1 \oplus \cdots \oplus A_n$
sont isomorphes. 
Pour ceci et d'autres exemples de d\'ecompositions non isomorphes,
voir les \S~90--91 de \cite{Fuch}~;
voir aussi \cite{Baum}, d\'ej\`a cit\'e dans l'introduction.

\medskip

Nous revenons \`a des exemples ind\'ecomposables,
et en particulier \`a certains groupes de Coxeter.

\medskip

\proclaim{Proposition 1}
Soit $\Gamma$ un groupe poss\'edant un sous-groupe normal ab\'elien $T$
avec les propri\'et\'es suivantes~:
\roster
\item"(i)"    il existe un entier $d \ge 1$ tel que
              $T$ est isomorphe \`a $\bold Z^d$~;
\item"(ii)"   l'action de $\Gamma/T$ sur $T$ est fid\`ele~;
\item"(iii)"  la repr\'esentation associ\'ee de $\Gamma/T$ sur
              $T \otimes_{\bold Z}\bold Q \approx \bold Q^d$ est irr\'eductible.
\endroster
Alors $\Gamma$ est ind\'ecomposable.
\endproclaim

\demo{D\'emonstration}
Montrons d'abord que tout sous-groupe normal ab\'elien $N$ de $\Gamma$ est
contenu dans $T$. 

Le sous-groupe $[N,T]$ de $\Gamma$ est dans $T$ et $\Gamma/T$-invariant.
La propri\'et\'e (iii) implique qu'il est ou bien d'indice fini dans $T$
ou bien r\'eduit \`a un \'el\'ement.
S'il \'etait d'indice fini, il existerait un entier $k \ge 1$ tel que
$[N,T]$ contienne un sous-groupe $kT \approx k\bold Z^d$~;
par la propri\'et\'e (ii), $N$ agirait non trivialement sur $[N,T]$
(qui est dans $N$), ce qui est impossible puisque $N$ est ab\'elien.
Donc $[N,T]=\{1\}$~; il en r\'esulte que l'action de $N$ sur $T$ est triviale,
de sorte que $N \subset T$ par la propri\'et\'e (ii).

Soient $\Gamma_1,\Gamma_2$ des sous-groupes de $\Gamma$ tels que 
$\Gamma = \Gamma_1 \times \Gamma_2$.
Posons
$$
\aligned
N_1 \, &= \, \{a \in \Gamma_1 \mid 
   \text{il existe} \quad b \in \Gamma_2 
   \quad \text{tel que} \quad (a,b) \in T\} , \\
N_2 \, &= \, \{b \in \Gamma_2 \mid 
   \text{il existe} \quad a \in \Gamma_1 
   \quad \text{tel que} \quad (a,b) \in T\} , \\
N \, &= \, N_1 \times N_2 .
\endaligned
$$
Alors $N_1$ est un sous-groupe normal ab\'elien de $\Gamma_1$
et $N_2$ un sous-groupe normal ab\'elien de $\Gamma_2$,
donc $N$ est un sous-groupe normal ab\'elien de $\Gamma$~;
de plus, $N$ contient $T$.
Il r\'esulte de la maximalit\'e de $T$ \'etablie plus haut que $N = T$.
La propri\'et\'e (iii) implique  que l'un des facteurs $N_1$, $N_2$
est r\'eduit \`a un \'el\'ement~; convenons que $N_2 = \{1\}$.
Comme $\Gamma_2$ centralise $N_1$, la propri\'et\'e (ii) implique que 
$\Gamma_2 = \{1\}$.
\hfill $\square$
\enddemo

\proclaim{Corollaire}
Un groupe de Coxeter irr\'eductible de type affine est ind\'ecomposable.
\endproclaim

\demo{D\'emonstration} 
Un groupe de Coxeter $W_a$ de type affine
est un produit semi-direct $Q \rtimes W$,
o\`u $W$ est un groupe de Weyl, en particulier un groupe fini,
et $Q$ le groupe des poids radiciels correspondant,
en particulier un groupe ab\'elien libre de type fini
(\cite{BouL}, chapitre~6, \S~2, no~1, proposition~1).
De plus, les conditions suivantes sont \'equivalentes~:
(i) $W_a$ est irr\'eductible comme groupe de Coxeter,
(ii) $W$ est irr\'eductible comme groupe de Coxeter,
(iii) la repr\'esentation de $W$ dans $Q \otimes_{\bold Z}\bold Q$ est
irr\'eductible (\cite{BouL}, chapitre~6, \S~1, no~1 et chapitre~5, \S~3, no~7).

La proposition 1 s'applique donc \`a la situation du corollaire.
\hfill $\square$
\enddemo

\bigskip
\head{\bf
3.~D\'ecompositions de sous-groupes de groupes alg\'ebriques
}\endhead

   Il est bien connu qu'un groupe dont le centre est r\'eduit \`a un \'el\'ement
poss\`ede au plus une d\'ecomposition en produit direct d'un nombre fini de
sous-groupes ind\'ecomposables. Nous commen\c cons par rappeler ce r\'esultat et
quelques-unes de ses cons\'equences, notamment le fait que la propri\'et\'e
d'ind\'ecomposabilit\'e convenablement formul\'ee passe de certains groupes
topologiques \`a leurs sous-groupes denses (voir aussi la proposition 9).

   Nous notons $Z(H)$ le centre d'un groupe $H$.
Si $H$ est sous-groupe d'un groupe $G$, 
nous \'ecrivons $H'$ son centralisateur dans $G$.
Lorsque $G$ est un produit direct $A \times B$, nous avons $A' = Z(A) \times B$.
Si de plus $Z(G) = \{1\}$, alors $A' = B$ et $A = B'$.

\medskip

\proclaim{Proposition 2} Soient $G$ un groupe de centre r\'eduit \`a un \'el\'ement
et $G_1,\hdots,G_n$ des sous-groupes ind\'ecomposables de $G$
tels que $G = G_1 \times \cdots \times G_n$.

Si $A,B$ sont deux sous-groupes de $G$ tels que $G = A \times B$,
il existe une renum\'erotation des $G_i$
et un entier $m \in \{0,\hdots,n\}$ tels que
$$ 
A \, = \, G_1 \times \cdots \times G_m
\qquad \text{et} \qquad
B \, = \, G_{m+1} \times \cdots \times G_n .
$$
\endproclaim

\demo{D\'emonstration}
Notons $g = (g_1,\hdots,g_n)$ l'\'ecriture d'un \'el\'ement $g \in G$
selon la d\'ecompo\-sition $G = G_1 \times \cdots \times G_n$.
Pour tout $i \in I \Doteq \{1,\hdots,n\}$, notons
$\pi_i : G \longrightarrow G_i , g \longmapsto g_i$
la projection canonique.
Posons $A_i = A \cap G_i$ et $B_i = B \cap G_i$.

Nous affirmons que $A = A_1 \times \cdots \times A_n$, de sorte que
$A_i = \pi_i(A)$ pour tout $i \in I$.
En effet, soit $a = (a_1,\hdots,a_n) \in A$.
Comme $a \in B'$, 
nous avons 
$$
[(a_1,\hdots,a_n),(b_1,\hdots,b_n)] \, = \, 
([a_1,b_1],\hdots,[a_n,b_n]) \, = \, (1,\hdots,1)
$$
pour tout $(b_1,\hdots,b_n) \in B$.
Il en r\'esulte que, pour $i \in I$, nous avons aussi
$[a_i,b_j] = 1$ pour tous $j \in J$ et $b_j \in B_j$,
et par suite $a_i \in B' = A$.
L'affirmation en r\'esulte.
De m\^eme $B = B_1 \times \cdots \times B_n$ 
et $B_i = \pi_i(B)$ pour tout $i \in I$.

Soit $i \in I$~; vu que $G = A \times B$ et $[A,B] = 1$, 
nous avons $G_i = A_i \times B_i$. 
De plus, l'un des groupes $A_i,B_i$ est r\'eduit \`a un \'el\'ement
parce que $G_i$ est ind\'ecomposable. La proposition en r\'esulte.
\hfill $\square$
\enddemo

\proclaim{Cons\'equence} 
    Un groupe $G$ qui satisfait aux hypoth\`eses de la proposition~2 
est bien s\^ur uniquement directement d\'ecomposable,
mais les sous-groupes ind\'ecomposables dont $G$ est le produit direct
sont de plus uniquement d\'etermin\'es comme sous-groupes 
(et non pas seulement \`a isomorphisme pr\`es).
\endproclaim

    L'argument de la preuve pr\'ec\'edente est suffisamment robuste pour
s'adapter \`a d'autres cas. 
Nous consid\'erons ci-dessous un corps alg\'ebriquement clos $K$
et des groupes alg\'ebriques d\'efinis sur $K$.
Un tel groupe $G$ est {\it ind\'ecomposable} s'il n'est pas produit direct 
{\it de sous-groupes alg\'ebriques} de mani\`ere non banale.
Pour des raisons de dimension, tout $G$ est produit direct d'une famille
finie de sous-groupes alg\'ebriques ind\'ecomposables~; 
lorsque de plus $Z(G) = \{1\}$, 
ces sous-groupes sont uniquement d\'etermin\'es \`a l'ordre pr\`es.

\medskip

\proclaim{Proposition 3}
Soient $G$ un groupe alg\'ebrique 
de centre r\'eduit \`a un \'el\'ement
et $G_1,\hdots,G_n$ 
des sous-groupes alg\'ebriques ind\'ecomposables de $G$
tels que $G = G_1 \times \cdots \times G_n$.
Soit $\Gamma$ un sous-groupe Zariski-dense de $G$.

Si $A,B$ sont deux sous-groupes de $\Gamma$ tels que $\Gamma = A \times B$,
il existe une renum\'erotation des $G_i$
et un entier $m \in \{0,\hdots,n\}$ tels que
$$ 
A \, \subset \, G_1 \times \cdots \times G_m
\qquad \text{et} \qquad
B \, \subset \, G_{m+1} \times \cdots \times G_n .
$$

En particulier, si $G$ est de plus ind\'ecomposable comme groupe alg\'ebrique,
c'est-\`a-dire si $n=1$, 
tout sous-groupe Zariski-dense de $G$ est ind\'ecomposable.
\endproclaim

\demo\nofrills{D\'emonstration (voir aussi l'exemple (x) du chapitre 2).
\hskip-1cm}  
Notons $\overline{A}$ et $\overline{B}$
les adh\'erences de Zariski de $A$ et $B$. 
Alors $[\overline{A},\overline{B}] = \{1\}$
(\cite{Bore}, no 2.4)~;
de plus, $\overline{A}  \hskip.1cm \overline{B}$ est ferm\'e dans $G$
(\cite{Bore}, no 1.4), de sorte que $G = \overline{A} \hskip.1cm \overline{B}$.
Comme $\overline{A} \cap \overline{B}$ est central dans $G$, 
cette intersection est r\'eduite \`a un \'el\'ement 
et $G$ est produit direct de $\overline{A}$ et $\overline{B}$.

L'argument de la d\'emonstration pr\'ec\'edente montre que, 
apr\`es renum\'erotation \'eventuelle des $G_j$, 
il existe un entier $m \in \{0,\hdots,n\}$ tel que
$A \subset \overline{A} = G_1 \times \cdots \times G_m$ et
$B \subset \overline{B} = G_{m+1} \times \cdots \times G_n$.
\hfill $\square$
\enddemo

\demo{Exemple}
Consid\'erons le groupe $\Gamma = \bold Z \ltimes \bold Z [1/pq]$,
o\`u  $p \ne 0$, $q \ge 2$ sont deux entiers premiers entre eux,
et o\`u le g\'en\'erateur $1$ de $\bold Z$ agit sur $\bold Z [1/pq]$
par multipication par~$p/q$.
C'est un groupe m\'etab\'elien de type fini~;
dans le cas o\`u $\abs{p} = 1$, 
c'est aussi un groupe de Baumslag-Solitar.
Notons $G$ le groupe des matrices triangulaires de la forme
$\left( \matrix a & b \\ 0 & 1 \endmatrix \right)$,
avec $a \in \bold C^*$ et $b \in \bold C$~;
c'est un groupe alg\'ebrique connexe ind\'ecomposable
(parce que son alg\`ebre de Lie l'est,
puisque c'est l'alg\`ebre de Lie r\'esoluble non ab\'elienne 
de dimension~$2$, voir ci-dessous avant la proposition 6).
Comme $\Gamma$ est Zariski-dense dans $G$, 
il r\'esulte de la proposition~3 que $\Gamma$
est ind\'ecomposable.
En particulier, tout groupe de Baumslag-Solitar r\'esoluble
qui n'est pas ab\'elien libre de rang deux  est ind\'ecomposable.
(Les autres groupes de Baumslag-Solitar sont \'egalement ind\'ecomposables~;
voir l'exemple (xii) du chapitre 2.)
\enddemo

\medskip

Avant de g\'en\'eraliser ces propositions \`a des cas avec centres,
nour rappelons les points suivants concernant la notion d'hypercentre.

Consid\'erons \`a nouveau un groupe $G$, sans autre structure.
La {\it suite centrale ascendante}  de $G$
est la suite $(\zeta^{\alpha}(G))_{\alpha}$,
index\'ee par les ordinaux $\alpha$,
d\'efinie par r\'ecurrence transfinie comme suit
(o\`u $\pi_{\beta}$ d\'esigne 
la projection canonique de $G$ sur $G/\zeta^{\beta}(G)$)~:

   si $\alpha = 0$, alors $\zeta^0(G) = \{1\}$~;

   si $\alpha = \beta+1$, alors
   $\zeta^{\alpha}(G) = \pi_{\beta}^{-1}\big(Z(G/\zeta^{\beta}(G))\big)$~;

   si $\alpha$ est un ordinal limite, alors
$\zeta^{\alpha}(G) = \bigcup_{\beta < \alpha} \zeta^{\beta}(G)$.

\noindent
Chaque $\zeta^{\alpha}(G)$ est un sous-groupe caract\'eristique de $G$.
{\it L'hypercentre} de $G$ est la r\'eunion $\zeta^{\uparrow}(G)$ des
$\zeta^{\alpha}(G)$~; 
il suffit de prendre la r\'eunion sur l'ensemble des ordinaux
dont le cardinal ne d\'epasse pas celui de $G$.
L'hypercentre de $G$ est un sous-groupe normal tel que
$Z(G/\zeta^{\uparrow}(G)) = \{1\}$, et qui est minimal pour cette propri\'et\'e.
En particulier, $\zeta^{\uparrow}(G) = \{1\}$ si et seulement si $Z(G) = \{1\}$.

Soit maintenant $G$ un groupe alg\'ebrique {\it connexe}.
Son centre $Z(G) = \zeta^1(G)$ est ou bien de dimension strictement positive 
ou bien fini. 
Dans le second cas, $\zeta^2(G)$ est ou bien de dimension strictement
positive ou bien fini, et alors \'egal \`a $Z(G)$ 	
puisque tout sous-groupe normal fini d'un groupe connexe est central.
Plus g\'en\'eralement, pour tout entier $k \ge 1$, l'une au moins des deux relations
$$
\dim\left( \zeta^k(G) \right) \, > \, \dim\left( \zeta^{k-1}(G) \right)
, \hskip1cm 
\zeta^{k+1} (G) \, = \, \zeta^k (G)
$$
est vraie. 
Il en r\'esulte qu'il existe un entier $h$ tel que 
$\zeta^{\uparrow}(G) = \zeta^h(G)$~; 
en particulier l'hypercentre de $G$ est un sous-groupe alg\'ebrique de $G$.
De m\^eme, l'hypercentre d'un groupe de Lie r\'eel ou complexe connexe
est un sous-groupe ferm\'e.

\medskip

\proclaim{Proposition 4} Soient $G$ un groupe, 
$H = G/\zeta^{\uparrow}(G)$ le quotient de $G$ par son hypercentre
et $\pi : G \longrightarrow H$ la projection canonique.
Soient $H_1,\hdots,H_n$ des sous-groupes ind\'ecomposables de $H$ tels que
$H = H_1 \times \cdots \times H_n$~;
posons $G_i = \pi^{-1}(H_i)$,
de sorte que $G = G_1 \cdots G_n$ et
$G_i \cap \prod_{j\ne i}G_j = \zeta^{\uparrow}(G)$ pour tout $i \in \{1,\hdots,n\}$.

Si $A,B$ sont deux sous-groupes de $G$ tels que $G = A \times B$,
il existe une renum\'erotation des $G_i$ 
et un entier $m \in \{0,\hdots,n\}$ tels que
$$ 
A \, \subset \, G_1 \cdots G_m
\qquad \text{et} \qquad
B \, \subset \, G_{m+1}  \cdots  G_n .
$$
\endproclaim

\demo{D\'emonstration} Il est \'evident que $H = \pi(A)\pi(B)$.
Par ailleurs, 
tout \'el\'ement de $\pi(A)$ commute \`a tout \'el\'ement de $\pi(B)$~;
par suite, tout \'el\'ement de $\pi(A) \cap \pi(B)$
commute \`a tout \'el\'ement de $\pi(A)\pi(B)$.
Comme $Z(H) = \{1\}$, il en r\'esulte que $H$ est produit direct
de ses sous-groupes $\pi(A)$ et $\pi(B)$.

La proposition 4 est donc une cons\'equence imm\'ediate de la proposition~2.
\hfill $\square$
\enddemo

\medskip

\proclaim{Proposition 5} Soient $G$ un groupe alg\'ebrique,
$H = G/\zeta^{\uparrow}(G)$ le quotient de $G$ par son hypercentre
et $\pi : G \longrightarrow H$ la projection canonique.
Soit $H = H_1 \times \cdots \times H_n$
la d\'ecomposition canonique de $H$ 
en produit de sous-groupes alg\'ebriques ind\'ecomposables~;
posons $G_i = \pi^{-1}(H_i)$, de sorte que
$G = G_1\cdots G_n$ et $G_i \cap \prod_{j \ne i}G_j = \zeta^{\uparrow}(G)$
pour tout $i \in \{1,\hdots,n\}$. 
Soit $\Gamma$ un sous-groupe Zariski-dense de $G$.

Si $A,B$ sont deux sous-groupes de $\Gamma$ tels que $\Gamma = A \times B$,
il existe une renum\'erotation des $G_i$ et un entier $m \in \{0,\hdots,n\}$
tels que
$$ 
A \, \subset \, G_1 \cdots G_m
\qquad \text{et} \qquad
B \, \subset \, G_{m+1}  \cdots  G_n .
$$

En particulier, si de plus $H$ est ind\'ecomposable comme groupe
alg\'ebrique, l'un des facteurs $A,B$ est contenu dans l'hypercentre
$\zeta^{\uparrow}(G)$.
\endproclaim

\demo{D\'emonstration}
Nous avons $\pi(\overline{A}) = \overline{\pi(A)}$,
$\pi(\overline{B}) = \overline{\pi(B)}$
et $[ \overline{\pi(A)}, \overline{\pi(B)}] = \{1\}$
(\cite{Bore}, num\'eros 1.4 et 2.4).
L'argument de la d\'emonstration de la proposition 3 montre que
$H = \pi(A) \times \pi(B)$.

La proposition 5 est alors une cons\'equence imm\'ediate de la proposition~3.
\hfill $\square$
\enddemo

\demo{Exemple} La proposition 5 s'applique \`a un groupe r\'eductif $G$
tel que $G/Z(G)$ soit simple, donc en particulier au groupe
$\operatorname{GL}_n(K)$~; c'est alors la proposition 10 de \cite{Mari}.
\enddemo

\demo{Exemple} Un groupe alg\'ebrique connexe ind\'ecomposable
dont le centre n'est pas r\'eduit \`a un \'el\'ement
peut contenir un sous-groupe Zariski-dense d\'ecomposable.
Consid\'erons en effet le groupe $G = \operatorname{SL}_4(\bold C)$. 
Notons $A$ le noyau de la r\'eduction
$\operatorname{SL}_4(\bold Z) \longrightarrow \operatorname{SL}_4(\bold Z/3\bold
Z)$  mo\-dulo~$3$
et $B$ le centre $\{\pm 1\}$ de $\operatorname{SL}_4(\bold Z)$. 
Alors le produit direct $\Gamma = A \times B$ 
est un sous-groupe Zariski-dense de $G$.
(Voir n\'eanmoins le chapitre 5 sur la c-ind\'ecomposabilit\'e.)
\enddemo

Pour appliquer les propositions 3 et 5, il convient de disposer du crit\`ere
d'ind\'ecompo\-sabilit\'e suivant pour les groupes alg\'ebriques. 
Rappelons qu'une alg\`ebre de Lie $\eufm{g}$ est {\it ind\'ecom\-posable} 
si $\eufm{g} \ne \{0\}$ et si,
pour tout isomorphisme de $\eufm{g}$ avec un produit d'alg\`ebres de Lie
$\eufm{a} \times \eufm{b}$, 
l'une de $\eufm{a}, \eufm{b}$ est r\'eduite \`a z\'ero.

\proclaim{Crit\`ere} Pour qu'un groupe alg\'ebrique connexe $G$ 
soit ind\'ecomposable,  il suffit que son alg\`ebre de Lie le soit.
\endproclaim

\proclaim{Strat\'egie} Soient $K$ un corps alg\'ebriquement clos,
$G$ un groupe alg\'ebrique connexe \`a centre fini
dont l'alg\`ebre de Lie $\eufm{g}$ est ind\'ecomposable,
$\Gamma$ un sous-groupe de $G$
et $\Gamma = A \times B$ une d\'ecomposition en produit direct.

Il r\'esulte du crit\`ere ci-dessus et des propositions qui pr\'ec\`edent que, 
si $\Gamma$ est Zariski-dense dans $G$, 
l'un des groupes $A,B$ est central d'ordre major\'e par l'ordre du centre de~$G$.

Supposons de plus que $G$ est un $L$-groupe,
o\`u $L$ est un sous-corps de $K$,
et que $\Gamma$ est un sous-groupe du groupe $G(L)$ des points rationnels.
Supposons aussi que $L$ est un corps parfait, 
de sorte que $G(L)$ est Zariski-dense dans $G$ (corollaire 18.3 de \cite{Bore}).
Si $\Gamma$ est Zariski-dense dans $G(L)$,
alors de m\^eme l'un des groupes $A,B$ est central d'ordre major\'e par
l'ordre du centre de $G(L)$.
\endproclaim

\medskip   

   Il y a des exemples imm\'ediats d'alg\`ebres de Lie ind\'ecomposables~:
une alg\`ebre de Lie simple, 
une alg\`ebre de Lie non ab\'elienne de dimension $2$, 
une alg\`ebre de Lie nilpotente \`a centre de dimension $1$.
Voici deux autres familles~: 
pour tout entier $n \ge 1$, le produit semi-direct standard
$\bold C^n \rtimes \eufm{gl}_n(\bold C)$ est ind\'ecomposable~;
une alg\`ebre de Lie isomorphe \`a une sous-alg\`ebre de Lie parabolique
d'une alg\`ebre de Lie complexe simple est ind\'ecomposable~;
voir les th\'eor\`emes 4.2 et 4.7 de \cite{Meng}.

   Nous verrons au num\'ero suivant une famille d'exemples d'alg\`ebres de Lie
r\'eelles ind\'e\-composables permettant d'appliquer aux groupes
de Coxeter la strat\'egie ci-dessus.

\medskip

\noindent{\it Digression.}
   Une alg\`ebre de Lie $\eufm{g}$ de dimension finie (sur un corps arbitraire)
s'\'ecrit comme produit direct d'alg\`ebres ind\'ecom\-posables
(la v\'erification par r\'ecurrence sur la dimension est de pure routine).
De plus, il y a unicit\'e au sens de la proposition suivante,
du type Wedderburn-Remak-Krull-Schmidt-Azumaya.
Bien que ce soit un r\'esultat classique,
nous n'avons pas su en trouver dans la litt\'erature la formulation qui nous
convient.

\medskip

\proclaim{Proposition 6} Soient $\eufm{g}$ une alg\`ebre de Lie de dimension
finie sur un corps $K$ et 
$$
\eufm{g} \, = \, \eufm{a}_1 \oplus \cdots \oplus \eufm{a}_m
\, = \, \eufm{b}_1 \oplus \cdots \oplus \eufm{b}_n
\tag*
$$
deux d\'ecompositions de $\eufm{g}$ 
en produits directs d'id\'eaux ind\'ecomposables. 

Alors $m=n$ et il existe une permutation $\sigma$ de $\{1,\hdots,m\}$ telle que
$\eufm{a}_j$ et $\eufm{b}_{\sigma(j)}$ sont isomorphes
pour tout $j \in \{1,\hdots,m\}$.
\endproclaim
\noindent

\demo{D\'emonstration}
C'est un r\'esultat tout \`a fait standard 
qu'il y a unicit\'e de la d\'ecomposition de $\eufm{g}$ en
{\it sous-$\eufm{g}$-modules} ind\'ecomposables
(voir par exemple le \lq\lq th\'eor\`eme
de Krull-Schmidt-Azumaya\rq\rq , no 19.21 dans \cite{Lam}, 
appliqu\'e \`a l'alg\`ebre enveloppante de $\eufm{g}$).

Pour une application lin\'eaire entre deux id\'eaux de $\eufm{g}$,
les conditions d'\^etre un morphisme de $\eufm{g}$-module
et d'\^etre un morphisme d'alg\`ebres de Lie sont {\it diff\'erentes}.
Toutefois, les isomorphismes 
fournis par {\it la preuve} du th\'eor\`eme invoqu\'e
sont des compositions d'injections et de projections canoniques
associ\'ees aux d\'ecompositions de (*)~;
ce sont donc {\it \`a la fois} des morphismes de $\eufm{g}$-modules
et des morphismes d'alg\`ebres de Lie.
Par suite, le th\'eor\`eme standard fournit bien
des isomorphismes d'alg\`ebres de Lie 
$\eufm{a}_j \longrightarrow \eufm{b}_{\sigma(j)}$.
\hfill $\square$
\enddemo

\bigskip
\head{\bf
4.~Groupes de Coxeter et groupes d'Artin
}\endhead

Soient $E$ un espace vectoriel r\'eel de dimension finie 
et $B$ une forme bilin\'eaire sym\'etrique sur $E$. 
Notons $r_B$ la dimension du noyau
$$
\operatorname{Ker}(B) \, = \,
\left\{ v \in E \mid B(v,w) = 0 \ \text{pour tout} \ w \in E \right\}
$$
de $B$~;
soit $p_B$ [respectivement $q_B$]
la dimension maximale d'un sous-espace $U$ de $E$
tel que $B(u,u) > 0$ [resp. $B(u,u) < 0$]
pour tout $u \in U$, $u \ne 0$. 
On sait que $p_B + q_B + r_B$ est la dimension de $E$, d\'esormais not\'ee $n_B$,
et nous appelons {\it signature} de $B$
le triplet $(p_B,q_B,r_B)$. 
Consid\'erons le groupe alg\'ebrique
$$
\operatorname{Of}(B) \, = \, \left\{g \in \operatorname{GL}(E) 
\ \Bigg\vert \
\aligned
B(gv,gw) &= B(v,w) \ \text{pour tout} \ v,w \in E \ \text{et} \\
gv &= v \ \text{pour tout} \ v \in \operatorname{Ker}(B)
\endaligned
\right\} ,
$$
qui est un produit semi-direct de la forme
$(\bold R^{p_B+q_B})^{r_B} \rtimes O(p_B,q_B)$.

Son alg\`ebre de Lie~$\eufm{of}(B)$
est isomorphe \`a l'alg\`ebre de Lie des ($n_B \times n_B$)-matrices 
d'\'ecriture par blocs
$$
\left( \matrix 
a & b & 0 \\ c & d & 0 \\ x & y & 0
\endmatrix \right)
$$
relativement \`a la d\'ecomposition $n_B = p_B + q_B + r_B$, avec
$$
\left( \matrix 
{}^t a & {}^t c  \\ {}^t b & {}^t d 
\endmatrix \right)
\left( \matrix
I_p & \phantom{-}0 \\ 0 & -I_q 
\endmatrix \right)
+
\left( \matrix
I_p & \phantom{-}0 \\ 0 & -I_q 
\endmatrix \right)
\left( \matrix 
a & b  \\ c & d 
\endmatrix \right)
\, = \,
0
$$
(les pr\'efixes ${}^t$ indiquent des transpositions,
et $I_p$ la ($p \times p$)-matrice unit\'e).

Lorsque $p_B + q_B \ge 3$, on v\'erifie que le centre de $\operatorname{Of}(B)$
est le groupe $\pm \operatorname{Id}_E$ d'odre $2$,
et qu'il co\"{\i}ncide avec son hypercentre.

\medskip

\proclaim{Proposition 7} Conservons les notations ci-dessus~;
supposons de plus 
$$
p_B + q_B \ \ge \ 2 
\qquad \text{et} \qquad 
(p_B,q_B,r_B) \ \ne \ (4,0,0), \ (2,2,0), \ (0,4,0) .
$$
Alors l'alg\`ebre de Lie $\eufm{of}(B)$ est ind\'ecomposable.
\endproclaim

\demo{Remarques}
Si $p_B+q_B \le 1$, l'alg\`ebre de Lie $\eufm{of}(B)$ est ab\'elienne.
Par ailleurs, les alg\`ebres de Lie
$\eufm{so}(4) = \eufm{so}(3) \times \eufm{so}(3)$
et 
$\eufm{so}(2,2) = \eufm{so}(2,1) \times \eufm{so}(2,1)
\approx \eufm{sl}_2(\bold R) \times \eufm{sl}_2(\bold R)$
sont d\'ecomposables.
Notons que
$\eufm{so}(3,1) = \eufm{so}(1,3) \approx \eufm{sl}_2(\bold C)$ est simple 
(il faut voir ici $\eufm{sl}_2(\bold C)$ 
comme une alg\`ebre de Lie {\it r\'eelle}), 
m\^eme si l'alg\`ebre de Lie complexifi\'ee ne l'est pas.
\enddemo

\demo{D\'emonstration} 
Si $r_B = 0$ et $(p_B,q_B) \ \ne \ (4,0), \ (2,2), \ (0,4)$, 
l'alg\`ebre de Lie $\eufm{of}(B) = \eufm{so}(p_B,q_B)$ est 
ou bien simple (si $p_B + q_B \ge 3$) ou bien de dimension un (si $p_B + q_B = 2$),
donc ind\'ecom\-posable dans tous ces cas. 
Supposons d\'esormais $r_B \ge 1$. 
  
Comme les signatures $(q_B,p_B,r_B)$ et $(p_B,q_B,r_B)$ donnent lieu \`a des
alg\`ebres de Lie isomorphes, nous pouvons aussi supposer $p_B \ge q_B$
sans restreindre la g\'en\'eralit\'e de ce qui suit.
Soient $\eufm{a}, \eufm{b}$
deux id\'eaux de $\eufm{of}(B)$ tels que
$\eufm{of}(B) = \eufm{a} \times \eufm{b}$~;
il s'agit de montrer que $\eufm{a} = 0$ ou $\eufm{b} = 0$.
Nous s\'eparons la discussion en trois cas.

\medskip

{\it (i) Cas o\`u $p_B+q_B \ge 3$ et $(p_B,q_B) \ne (4,0), (2,2)$.}
Soit $\eufm{s}$ une sous-alg\`ebre de Levi de $\eufm{of}(B)$~;
c'est une alg\`ebre de Lie simple.
Sa projection sur $\eufm{a}$ est ou bien isomorphe \`a $\eufm{s}$
ou bien nulle, et de m\^eme pour sa projection sur $\eufm{b}$. 
Supposons les notations telles que sa projection sur $\eufm{b}$ soit nulle,
c'est-\`a-dire telles que l'alg\`ebre de Lie $\eufm{b}$ soit r\'esoluble.
Une v\'erification de routine \'etablit que 
l'alg\`ebre de Lie $\eufm{of}(B)$ est parfaite
(rappelons qu'une alg\`ebre de Lie est dite {\it parfaite}
si elle co\"{\i}ncide avec son id\'eal d\'eriv\'e).
L'alg\`ebre de Lie $\eufm{b}$ est donc r\'esoluble et parfaite,
c'est-\`a-dire r\'eduite \`a z\'ero.

\medskip

{\it (ii) Cas o\`u $p_B+q_B = 2$.}
L'alg\`ebre $\eufm{of}(B)$, r\'esoluble, correspond aux matrices de la forme
$
\left( \matrix 
0 & c & 0 \\ c & 0 & 0 \\ x & y & 0
\endmatrix \right)
$
avec $c \in \bold R$ et $x,y \in \bold R^{r_B}$.
Son radical nilpotent $\eufm{n}$, 
correspondant aux matrices pour lesquelles $c=0$,
est de codimension $1$~;
il se projette donc surjectivement sur au moins l'un des facteurs
$\eufm{a}, \eufm{b}$.
Ce facteur \'etant \`a la fois nilpotent et parfait 
(puisque $\eufm{n} = [\eufm{of}(B),\eufm{of}(B)]$),
il est r\'eduit \`a z\'ero,
et $\eufm{of}(B)$ est bien ind\'ecomposable.

\medskip

{\it (iii) Cas o\`u  $(p_B,q_B) = (4,0)$ ou $(p_B,q_B) = (2,2)$.}
Notons \`a nouveau $\eufm{s}$ une sous-alg\`ebre de Levi de $\eufm{of}(B)$~;
elle poss\`ede deux id\'eaux simples isomorphes $\eufm{u}$ et $\eufm{v}$
tels que $\eufm{s} = \eufm{u} \times \eufm{v}$.
Vu l'argument du cas (i), il suffit de consid\'erer ici le cas o\`u
la projection de $\eufm{s}$ sur $\eufm{a}$ 
serait isomorphe \`a $\eufm{u}$
et sa projection sur $\eufm{b}$ isomorphe \`a $\eufm{v}$~;
nous allons montrer que ce cas ne se produit pas.

   En effet, le facteur $\eufm{a}$ poss\`ederait 
une sous-alg\`ebre simple de dimension trois 
centra\-lisant le radical r\'esoluble du facteur $\eufm{b}$,
et de m\^eme pour le facteur $\eufm{b}$ 
et sa sous-alg\`ebre simple~$\eufm{v}$
centralisant le radical r\'esoluble de $\eufm{a}$.
La repr\'esentation naturelle (par restriction de la repr\'esentation adjointe) de
$\eufm{s}$  sur le radical r\'esoluble de $\eufm{of}(B)$
contiendrait donc des sous-repr\'esentations irr\'eductibles non fid\`eles.
Or ceci est absurde, car les sous-espaces irr\'eductibles 
de l'action de $\eufm{s} = \eufm{so}(p_B,q_B)$
sur le radical r\'esoluble $(\bold R^{p_B+q_B})^{r_B}$
sont tous isomorphes au 
$\eufm{so}(p_B,q_B)$-module fid\`ele $\bold R^{p_B+q_B}$.
\hfill $\square$
\enddemo

\medskip

Soit $(W,S)$ un syst\`eme de Coxeter, avec $S$ fini.
Notons $E$ l'espace vectoriel $\bold R^S$ 
et $B$ la forme de Tits associ\'ee \`a $(W,S)$.
Alors $W$ poss\`ede une {\it repr\'esentation g\'eom\'etrique} sur $E$ 
pour laquelle la forme $B$ est invariante.
Cette repr\'esentation est fid\`ele, 
et fournit donc une injection $W \subset \operatorname{Of}(B)$,
o\`u $\operatorname{Of}(B)$ est le groupe alg\'ebrique introduit plus haut~;
de plus $W$ est un sous-groupe discret 
du groupe des points r\'eels de $\operatorname{Of}(B)$.

Le groupe $W$ est dit {\it irr\'eductible} 
si le graphe de Coxeter associ\'e au syst\`eme $(W,S)$ est connexe.
Dans ce cas, la repr\'esentation g\'eom\'etrique de $W$ est ind\'ecomposable~;
de plus, les trois conditions suivantes sont \'equivalentes~:
cette repr\'esentation est irr\'eductible,
elle est absolument irr\'eductible,
le noyau de $B$ est r\'eduit \`a z\'ero.

Pour tout ceci, voir \cite{BouL}, chapitre 5, \S \ 4.

Rappelons quelques propri\'et\'es de la signature $(p_B,q_B,r_B)$ de $B$
lorsque $W$ est irr\'e\-ductible~; 
$n_B = p_B+q_B+r_B$ d\'esigne comme plus haut la dimension de $E$,
c'est-\`a-dire le cardinal de $S$.
\roster
\item"$\circ$" $p_B = n_B$ si et seulement si $W$ est fini~;
\item"$\circ$" $p_B = n_B - 1$ et $r_B = 1$ si et seulement si 
   $W$ est infini et contient un sous-groupe
   ab\'elien libre d'indice fini
   ($W$ est alors dit {\it de type affine})~;
\item"$\circ$" si $q_B = 0$, alors $r_B \le 1$~;
\item"$\circ$" si $n_B \le 4$, alors $p_B \ge n_B - 1$~;
   en particulier, $(p_B,q_B,r_B) \ne (2,2,0)$~;
\item"$\circ$" si $n_B \ge 4$, alors $p_B \ge 3$~;
\endroster
voir \cite{BouL} pour les trois premi\`eres propri\'et\'es, 
et 
\cite{Par1} pour les deux derni\`eres.

La proposition suivante appara\^{\i}t dans \cite{Par2}.

\medskip

\proclaim{Proposition 8} Soit $(W,S)$ un syst\`eme de Coxeter, avec $S$ fini.

(i) Si $W$ est irr\'eductible infini, 
$W$ est ind\'ecomposable.

(ii) Si $W$ est irr\'eductible infini non affine, 
tout sous-groupe d'indice fini de $W$ est ind\'ecomposable.

(iii) Dans tous les cas, $W$ est uniquement directement d\'ecomposable.
\endproclaim

\demo{D\'emonstration}
Notons que l'assertion (i) pour $W$ de type affine est une r\'ep\'etition du
corollaire de la proposition~1.

Soit $\Cal G$ le graphe de Coxeter associ\'e \`a la paire $(W,S)$.
Notons $\Cal G_1, \hdots, \Cal G_k$ les composantes connexes de ce graphe
et $W_1, \hdots, W_k$ les groupes de Coxeter correspondants,
qui sont les {\it composantes irr\'eductibles} du groupe $W$,
et dont $W$ est produit direct.
Ceux des $W_i$ qui sont infinis 
ont un centre r\'eduit \`a un \'el\'ement
({\cite{BouL}, chapitre 5, \S~4, exercice~3).
Si $W_i$ est infini et de plus n'est pas de type affine,
alors $W_i$ est Zariski-dense 
dans un groupe du type $\operatorname{Of}(B_i)$ \cite{BeHa}~;
il en r\'esulte en particulier que
le centre de tout sous-groupe d'indice fini de $W_i$
est encore r\'eduit \`a un \'el\'ement.

Il suffit donc d'appliquer la proposition 5 
pour d\'emontrer les assertions (i) et (ii).
L'assertion (iii) r\'esulte de la proposition 2
lorsque les $W_i$  sont des groupes de Coxeter \`a centres triviaux 
(par exemple sont tous des groupes de Coxeter infinis)~;
pour le cas g\'en\'eral,
nous invoquons la proposition 9. 
\hfill $\square$
\enddemo

La proposition suivante est un cas particulier du 
\lq\lq th\'eor\`eme fondamental\rq\rq \ du \S~47 de~\cite{Kuro}.

\medskip

\proclaim{Proposition 9}
Consid\'erons un entier $m \ge 1$, 
des groupes ind\'ecom\-posables $\Gamma_1,\hdots,\Gamma_m$
et le produit direct $\Gamma = \Gamma_1 \times \cdots \times \Gamma_m$~;
supposons 
\footnote{
Cette hypoth\`ese ne peut en aucun cas \^etre omise,
comme le montrent plusieurs des exemples d\'ej\`a cit\'es
(nos (xiv) et (xvii) du chapitre 2).
Le point important est que tout sous-groupe du centre de $\Gamma$
poss\`ede une s\'erie principale de sous-groupes normaux.
}
que le centre de $\Gamma$ est fini.
Alors $\Gamma$ est uniquement directement d\'ecomposable.
\endproclaim

Soient $W$ un groupe de Coxeter fini irr\'eductible d'un des types $A, D, E$
et $\Gamma$ le groupe d'Artin correspondant, de quotient $W$. 
I.~Marin a montr\'e qu'il existe 
un entier $d$ et une repr\'esentation irr\'eductible
de $\Gamma$ dans $\operatorname{GL}_d(\bold C)$ d'image Zariski-dense.
La proposition~5 permet donc une autre d\'emonstration du r\'esultat 
suivant de~\cite{Mari}.

\medskip

\proclaim{Proposition 10}
Soit $\Gamma_0$ un sous-groupe d'indice fini dans un groupe d'Artin $\Gamma$
comme ci-dessus. Si $\Gamma_0$ est directement d\'ecomposable, 
alors $\Gamma_0$ est produit direct de deux facteurs ind\'ecomposables exactement, 
dont l'un est central.
\endproclaim

   En utilisant un autre type d'argument, Paris a montr\'e que
tout groupe d'Artin irr\'educ\-tible de type sph\'erique est ind\'ecomposable 
(voir la proposition 4.2 de \cite{Par3}).

\medskip

  Consid\'erons en particulier le cas 
du groupe d'Artin $B_n$ des tresses \`a $n \ge 2$ brins.
On sait que le centre de $B_n$ est cyclique infini~;
plus g\'en\'eralement, tout sous-groupe d'indice fini $\Gamma$ de $B_n$
poss\`ede un centre $Z(\Gamma) = Z(B_n) \cap \Gamma$ 
qui est d'indice fini dans $Z(B_n) \approx \bold Z$.
Par suite, si un tel groupe $\Gamma$ n'est pas ind\'ecomposable, 
alors $\Gamma$ est produit direct d'un groupe ind\'ecomposable
$\Gamma_1$  et de son centre isomorphe \`a $\bold Z$~;
de plus, $\Gamma_1$ est isomorphe au quotient de $\Gamma$ par son centre,
et sa classe d'isomorphisme est donc d\'etermin\'ee par celle de $\Gamma$.
En particulier, $\Gamma$ est uniquement directement d\'ecomposable.

   Ceci s'applique au groupe $P_n$ des tresses pures pour tout $n \ge 3$. 
En effet,  nous avons une suite d'extensions scind\'ees 
$$                   
\{1\}  \longrightarrow F_{n-1} \longrightarrow P_n 
\overset{\pi_n}\to{\longrightarrow} P_{n-1} 
\longrightarrow \{1\}
\tag$\sharp$
$$
telles que l'image par $\pi_n$ du centre de $P_n$ 
co\"{\i}ncide avec le centre de $P_{n-1}$
(le noyau $F_{n-1}$ de $\pi_n$ est un groupe
libre \`a $n-1$ g\'en\'erateurs).
D\'efinissons par r\'ecurrence une suite $(Q_n)_{n \ge 2}$
de sous-groupes des $P_n$
en posant $Q_2 = \{1\}$
et $Q_n = \pi_n^{-1}(Q_{n-1})$ pour $n \ge 3$.
(Notons que $Q_3 = P_3$ est un groupe non ab\'elien libre \`a deux g\'en\'erateurs.
Notons aussi que, pour $n \ge 3$, $Q_n$ d\'epend du brin choisi pour d\'efinir
l'extension $(\sharp)$ ci-dessus, et n'est donc pas uniquement d\'efini comme
sous-groupe de $P_n$.) 
Il est facile de v\'erifier que $P_n$ est produit direct
de son centre, cyclique infini, et du groupe $Q_n$,
ind\'ecomposable.

\bigskip
\head{\bf
5.~Variation sur la notion d'ind\'ecomposabilit\'e
}\endhead

Un groupe $\Gamma$ est dit {\it c-d\'ecomposable}
s'il existe deux sous-groupes normaux infinis $\Gamma_1,\Gamma_2$ de $\Gamma$
tels que $\Gamma_1 \cap \Gamma_2$ est fini 
et $\Gamma_1\Gamma_2$ d'indice fini dans $\Gamma$.
Un groupe qui n'est pas c-d\'ecomposable est dit {\it c-ind\'ecomposable}.
(Voir \cite{Marg}, chap. IX, no 2.2~; 
nous \'ecrivons \lq\lq c-d\'ecomposable\rq\rq \ 
o\`u Margulis \'ecrit \lq\lq d\'ecom\-posable\rq\rq .
La lettre \lq\lq c\rq\rq \ est l'initiale de \lq\lq commensurable\rq\rq .)

Soient $\Gamma$ un groupe, $\Gamma_0$ un sous-groupe d'indice fini,
et $F$ un sous-groupe normal fini de $\Gamma$. 
Il est facile de v\'erifier que les trois conditions suivantes sont \'equivalentes~:
$\Gamma$ est c-ind\'ecomposable, 
$\Gamma_0$ est c-ind\'ecomposable,
$\Gamma/F$ est c-ind\'ecomposable.

Soit $A$ un ensemble fini non vide. 
Pour tout $\alpha \in A$, 
soient $k_{\alpha}$ un corps local 
et $G_{\alpha}$ un $k_{\alpha}$-groupe alg\'ebrique connexe,
presque $k_{\alpha}$-simple,  
et tel que le groupe $G_{\alpha}(k_{\alpha})$ n'est pas compact.
Soient $G$ le produit direct $\prod_{\alpha \in A} G_{\alpha}(k_{\alpha})$
et $\Gamma$ un r\'eseau dans $G$.
Alors $\Gamma$ est c-ind\'ecomposable (comme groupe abstrait)
si et seulement si $\Gamma$ est irr\'eductible (comme r\'eseau dans $G$)~;
voir \cite{Marg}, chap. IX, no 2.3.

\medskip
\noindent
{\it Remarque.} Soit $\Delta$ un groupe r\'esiduellement fini 
qui poss\`ede un sous-groupe d'indice fini $\Gamma$
h\'er\'editairement ind\'ecomposable, 
ce qui veut dire que tout sous-groupe d'indice fini de $\Gamma$
est ind\'ecomposable.  
Alors $\Delta$ est c-ind\'ecomposable.

En effet, soient $\Delta_1,\Delta_2$ deux sous-groupes normaux de $\Delta$
tels que $\Delta_1 \cap \Delta_2$ est fini 
et $\Delta_1\Delta_2$ d'indice fini dans $\Delta$.
Il s'agit de montrer que l'un des groupes $\Delta_1,\Delta_2$ est fini.
Pour $j=1,2$, posons $\Gamma_j' = \Delta_j \cap \Gamma$.
Vu l'hypoth\`ese de finitude r\'esiduelle, 
nous pouvons choisir  
un sous-groupe d'indice fini $\Gamma_j$ de $\Gamma_j'$ ($j=1,2$)
de telle sorte que que $\Gamma_1 \cap \Gamma_2 = \{1\}$.
Quitte \`a remplacer \`a nouveau $\Gamma_j$ par un sous-groupe d'indice fini,
nous pouvons supposer de plus que $\Gamma_1 \times \Gamma_2$
est un sous-groupe de~$\Gamma$.
Les hypoth\`eses impliquent alors que l'un des groupes $\Gamma_1,\Gamma_2$
est fini, de sorte que le groupe $\Delta_j$ correspondant est aussi fini.

\bigskip
\head{\bf
6.~Conditions Max et Min
}\endhead

  Un groupe satisfait \`a la {\it condition Max-n} 
si toute cha\^{\i}ne ascendante de sous-groupes normaux de $\Gamma$ 
est ultimement stationnaire,
et \`a la condition {\it Min-n}
si toute cha\^{\i}ne des\-cendante de sous-groupes normaux
est ultimement stationnaire.
On d\'efinit de m\^eme les {\it conditions Max-fd et Min-fd},
en termes de facteurs directs.
Un groupe satisfaisant \`a l'une de ces quatre conditions
est \'evidemment produit direct d'un nombre fini de groupes ind\'ecomposables.

   Un groupe satisfaisant 
\footnote{
Autrement dit~: un groupe poss\'edant une suite de composition distingu\'ee
principale, selon la terminologie de Bourbaki
(\cite{BouA}, chapitre I, \S~4, exercice~17).
} 
aux deux conditions Max-n et Min-n
est uniquement directement d\'ecomposable~:
c'est le th\'eor\`eme de Wedderburn-Remak-Krull-Schmidt~; 
voir par exemple 
le dernier th\'eor\`eme du \S~47 de \cite{Kuro},
ou le th\'eor\`eme 6.36 de \cite{Rotm},
ou le th\'eor\`eme 4.8 de \cite{Suzu, Chap.~2}.  
Les exemples de Baumslag cit\'es dans l'introduction 
montrent qu'un groupe satisfaisant la seule condition Max-n 
n'est pas n\'ecessairement uniquement directement d\'ecomposable, 
puisqu'un groupe polycyclique satisfait m\^eme 
la condition de cha\^{\i}ne ascendante  
pour les sous-groupes (non n\'ecessairement normaux). 

\medskip

Notons qu'un groupe de Coxeter infini ne poss\`ede jamais la propri\'et\'e Min-n.
Plus g\'en\'eralement, 
un groupe infini r\'esiduellement fini ne poss\`ede pas cette propri\'et\'e.
Or les groupes de Coxeter sont r\'esiduellement finis~: 
c'est en effet un fait g\'en\'eral, 
connu sous le nom de \lq\lq lemme de Mal'cev\rq\rq ,
que tout groupe lin\'eaire de type fini est r\'esiduellement fini~;
pour l'esquisse d'un argument valant pour les groupes de Coxeter, voir
\cite{BouL}, chapitre 5, \S~4, exercice 9.

Nous allons montrer 
qu'un groupe de Coxeter qui n'est pas virtuellement ab\'elien
ne poss\`ede jamais la propri\'et\'e Max-n. 
Il r\'esulte n\'eanmoins de la proposition~8 qu'un groupe de Coxeter
poss\`ede les propri\'et\'es Min-fd et Max-fd.

\medskip

\proclaim{Lemme} 
(i) Un quotient d'un groupe qui satisfait la condition 
Max-n la satisfait aussi.

(ii) Un sous-groupe d'indice fini d'un groupe qui satisfait la condition
Max-n la satisfait aussi.

(iii) Un groupe libre non ab\'elien ne satisfait pas la condition Max-n.
\endproclaim

\demo{D\'emonstration}
L'assertion (i) est banale et l'assertion (ii) est un r\'esultat de \cite{Wil1}.

C'est une cons\'equence imm\'ediate de la d\'efinition 
qu'un groupe $\Gamma$ satisfait la condition Max-n 
si et seulement si tout sous-groupe normal de $\Gamma$ 
peut \^etre engendr\'e {\it comme sous-groupe normal} par un ensemble fini. 
Il r\'esulte de l'existence de groupes \`a deux g\'en\'erateurs 
qui ne sont pas de pr\'esentation finie \cite{Neum} 
que le groupe libre de rang deux ne satisfait pas la condition Max-n, 
et donc de (i) qu'un groupe libre non ab\'elien de rang quelconque
ne la satisfait pas non plus.
\hfill $\square$
\enddemo

\medskip

\proclaim{Proposition 11} 
Un groupe de Coxeter qui n'est pas virtuellement ab\'elien 
ne satisfait pas la condition Max-n.
\endproclaim

\demo{D\'emonstration}
Selon un r\'esultat \'etabli ind\'ependamment  dans \cite{Gonc} et \cite{MaVi}
un groupe de Coxeter qui n'est pas virtuellement ab\'elien 
poss\`ede un sous-groupe d'indice fini 
qui se surjecte sur un groupe libre non ab\'elien.
La proposition est alors une cons\'equence imm\'ediate du lemme pr\'ec\'edent.
\hfill $\square$
\enddemo

\bigskip
\head{\bf
7.~Groupes \`a quotients major\'es
}\endhead

Soit $n$ un entier, $n \ge 2$.
Un groupe $\Gamma$ est dit {\it \`a quotients $n$-major\'es}, 
ou {\it $n$-QM},
si tout sous-groupe de type fini $\Delta \ne \{1\}$ de $\Gamma$ 
poss\`ede un sous-groupe normal propre d'indice au plus $n$.
Un groupe $\Gamma$ est dit {\it \`a quotients major\'es}
s'il existe un entier $n \ge 2$ pour lequel il est \`a quotients $n$-major\'es.
Nous collectons quelques exemples et propri\'et\'es simples 
relatifs \`a cette notion.

\medskip

(i) Un groupe fini est \'evidemment un groupe \`a quotients major\'es.

Pour un nombre premier $p$, 
un groupe qui est r\'esiduellement un $p$-groupe fini est $p$-QM.
En particulier, 
les groupes ab\'eliens libres et les groupes non ab\'eliens libres sont $2$-QM.

Le groupe $\oplus_{p \in \bold P}\bold Z /p\bold Z$, o\`u $\bold P$ d\'esigne 
l'ensemble des nombres premiers, n'est pas QM. 
Plus g\'en\'eralement, un groupe contenant des sous-groupes simples d'ordres
arbitrairement grands n'est pas QM.

\medskip

(ii) Pour la d\'efinition de la propri\'et\'e QM, 
l'exemple qui suit montre qu'on ne pourrait pas
omettre la condition sur $\Delta$ d'\^etre de type fini
sans changer la notion.

Rappelons d'abord un fait \'el\'ementaire~:
soient $p$ un nombre premier et $A$ un groupe ab\'elien
dans lequel tout \'el\'ement d'une part est divisible par $p$
et d'autre part a un ordre qui est une puissance de $p$~;
alors $A$ est un groupe divisible.
En effet, soient $a \in A$ et $\ell$ un nombre premier distinct de $p$~;
soit $n$ tel que $a$ soit d'ordre $p^n$, et soient $s,t \in \bold Z$
tels que $s\ell + tp^n = 1$~; si on pose $b = a^s$,
alors $b^{\ell} = a^{s\ell}\left(a^{p^n}\right)^t = a$.

En particulier, le groupe $\Gamma = \bold Z[1/p] / \bold Z$ est divisible,
et n'a donc aucun sous-groupe propre d'indice fini.
Toutefois, un sous-groupe propre de type fini de $\Gamma$ 
est un sous-groupe $(p^{-k}\bold Z) / \bold Z$, 
qui est cyclique d'ordre $p^k$,
pour un entier $k \ge 1$ convenable.
Par suite, $\Gamma$ est $p$-QM.

\medskip

(iii) Il est \'evident que tout sous-groupe d'un groupe $n$-QM est aussi $n$-QM.
Une somme restreinte (finie ou infinie) de groupes est $n$-QM si et seulement si
chaque facteur est $n$-QM.

\medskip

(iv) Soit $\Gamma$ un groupe qui s'ins\`ere dans  une extension 
$$
\{1\} \quad \longrightarrow \quad
\Gamma '  \quad \longrightarrow \quad
\Gamma  \quad \overset{\pi}\to{\longrightarrow} \quad
\Gamma ''  \quad \longrightarrow \quad
\{1\} .
$$
Si $\Gamma '$ est $n'$-QM et si $\Gamma ''$ est $n''$-QM,
alors $\Gamma$ est $\max(n',n'')$-QM.

En effet, soit $\Delta$ un sous-groupe de type fini de $\Gamma$.
Si $\Delta \subset \Gamma'$, alors $\Delta$ poss\`ede un sous-groupe normal propre
d'indice fini au plus $n'$.
Sinon, soit $\Delta_0''$ 
un sous-groupe normal propre d'indice $k \le n''$
de $\pi(\Delta)$~;
alors $\pi^{-1}(\Delta_0'')$
est un sous-groupe normal propre d'indice $k$ dans $\Delta$.

En particulier, si un groupe $\Gamma$ 
poss\`ede un sous-groupe normal d'indice fini qui est \`a quotients major\'es,
alors $\Gamma$ est \'egalement \`a quotients major\'es.

\medskip

(v) Pour les groupes de type fini,
aucune des propri\'et\'es \lq\lq  QM\rq\rq \ et
\lq\lq r\'esiduellement fini\rq\rq \ n'implique l'autre,
comme le montrent les consid\'erations suivantes.

D'une part, soit $S$ un groupe non ab\'elien et $T$ un groupe infini.
Le produit en couronne $\Gamma = S \wr T$
n'est pas r\'esiduellement fini (th\'eor\`eme 3.2 de  \cite{Grue}).
Si $S$ et $T$ sont deux groupes 
qui sont de type fini et $n$-QM pour un entier $n$,
il r\'esulte de (iii) et (iv) que le groupe de type fini
$\Gamma$ est aussi $n$-QM.
[Voir aussi (vi) ci-dessous.]

D'autre part, 
tout groupe d\'enombrable r\'esiduellement fini
se plonge dans un groupe de type fini et r\'esiduellement fini \cite{Wil2}.
Il r\'esulte donc de (i) et (iii) qu'il existe 
des groupes de type fini r\'esiduellement finis
qui ne sont pas QM.

\medskip

(vi) Rappelons un exemple de groupe r\'esoluble de type fini  
d\^u \`a P.~Hall \cite{HalP}.

Soient $R$ est un anneau commutatif avec unit\'e,
dont nous notons $R^*$ le groupe des unit\'es,
et $R_0$ un sous-groupe additif de $R$.
Posons
$$
G(R) \, = \, 
\left( \matrix 1 & R & R \\ 0 & R^* & R \\ 0 & 0 & 1 \endmatrix \right)
\qquad \text{et} \qquad 
Z(R_0) \, = \, 
\left( \matrix 1 & 0 & R_0 \\ 0 & 1 & 0 \\ 0 & 0 & 1 \endmatrix \right) .
$$
Le groupe $G(R)$ est r\'esoluble de classe $3$,
son centre s'identifie \`a $Z(R)$,
et $Z(R_0)$ est un sous-groupe central de $G(R)$.

Soit $p$ un nombre premier~; le groupe multiplicatif $p^{\bold Z}$
est un sous-groupe d'indice deux dans le groupe des unit\'es de $\bold Z [1/p]$.
Le groupe
$$
G_+(\bold Z[1/p]) \, = \,
\left( \matrix 
1 & \bold Z [1/p] & \bold Z [1/p] \\ 
0 & p^{\bold Z} & \bold Z [1/p] \\
0 & 0 & 1 
\endmatrix\right) ,
$$
qui est d'indice deux dans $G(\bold Z [1/p])$, est engendr\'e par les trois matrices
$$
\left( \matrix
1 & 0 & 0 \\ 0 & p & 0 \\ 0 & 0 & 1
\endmatrix \right) , \quad
\left( \matrix
1 & 1 & 0 \\ 0 & 1 & 0 \\ 0 & 0 & 1
\endmatrix \right)  \quad \text{et} \quad
\left( \matrix
1 & 0 & 0 \\ 0 & 1 & 1 \\ 0 & 0 & 1
\endmatrix \right) .
$$
Notons $\tilde\alpha$ l'automorphisme ext\'erieur de $G_+(\bold Z[1/p])$
obtenu en conjugant les matrices 
par la matrice diagonale de coefficients diagonaux $p, 1, 1$.
On v\'erifie que $\tilde\alpha(Z(\bold Z)) = Z(p\bold Z)$
est d'indice $p$ dans $Z(\bold Z)$.

L'exemple de Hall est le quotient $H$ de $G_+(\bold Z[1/p])$ par $Z(\bold Z)$.
L'automorphisme $\tilde\alpha$ de $G_+(\bold Z[1/p])$
induit un endomorphisme $\alpha$ de $H$ qui est surjectif
de noyau $Z(p^{-1}\bold Z)/Z(\bold Z)$,
c'est-\`a-dire de noyau cyclique d'ordre $p$.
En particulier, le groupe $H$ est 
de type fini, 
r\'esoluble (c'est m\^eme une extension centrale d'un groupe m\'etab\'elien)
et non Hopfien.
Nous avons une suite exacte courte
$$
\{1\} 
\, \longrightarrow \,
\bold Z [1/p] / \bold Z
\, \longrightarrow \,
H
\, \longrightarrow \,
G_+(\bold Z [1/p]) / Z (\bold Z [1/p])
\, \longrightarrow \, 
\{1\} 
$$
dont le noyau est le groupe $p$-QM de l'exemple (ii).
Le quotient de la suite exacte, 
isomorphe \`a $\bold Z \ltimes (\bold Z[1/p])^2$,
est un groupe lin\'eaire de type fini,
et c'est donc aussi un groupe QM (argument direct, ou proposition 12 ci-dessous).
Il r\'esulte de (iv) que le groupe de type fini non Hopfien $H$ est un groupe QM.

\medskip

(vii) Notre int\'er\^et pour la propri\'et\'e QM vient du fait 
que les groupes lin\'eaires de type fini l'ont. 

\medskip

\proclaim{Proposition 12} Soit $\Gamma$ un groupe de type fini qui est lin\'eaire,
c'est-\`a-dire qui est un sous-groupe de $\operatorname{GL}_d(K)$
pour un entier $d$ et un corps $K$ convenables.
Alors il existe un nombre premier $p$ tel que 
$\Gamma$ poss\`ede un sous-groupe d'indice fini 
qui est r\'esiduellement un $p$-groupe fini.

En particulier, $\Gamma$ est un groupe \`a quotients major\'es.
\endproclaim

\demo{D\'emonstration}
Soit $A$ le sous-anneau de $K$ engendr\'e par les coefficients matriciels
des \'el\'ements d'un syst\`eme fini de g\'en\'erateurs de $\Gamma$~;
c'est un anneau commutatif int\`egre de type fini.
Soit $\eufm{m}$ un id\'eal maximal de $A$~;
le quotient $A/\eufm{m}$ est un corps fini
(voir par exemple \cite{BoAC'}, chapitre 5, \S \ 3, no 4, corollaire 1
du th\'eor\`eme 3)  
dont nous notons $p$ la caract\'eristique.
Pour tout entier $k \ge 0$, notons $N_k$ 
le noyau de l'application naturelle 
$\operatorname{GL}_d(A) \longrightarrow \operatorname{GL}_d(A/\eufm{m}^k)$.

L'intersection des id\'eaux $\eufm{m}^k$ est r\'eduite \`a z\'ero
(r\'esultat de Krull, 
voir par exemple \cite{BoAC}, chap. 3, \S~3, no 2),
et donc  l'intersection des sous-groupes $N_k$ est r\'eduite \`a un \'el\'ement.
Par ailleurs, le quotient $N_k/N_{k+1}$ est fini pour tout $k \ge 0$,
et c'est un $p$-groupe ab\'elien \'el\'ementaire pour tout $k \ge 1$.
Il en r\'esulte que $N_1$ est r\'esiduellement  un $p$-groupe fini
qui est d'indice fini dans $\operatorname{GL}_d(A)$,
et par cons\'equent que $\Gamma \cap N_1$ est de m\^eme
r\'esiduellement un $p$-groupe fini qui est d'indice fini dans $\Gamma$. 

La seconde assertion de la proposition r\'esulte alors des points (i) et (iv)
ci-dessus.
\hfill $\square$
\enddemo

\medskip
\head{\bf
8.~Majorations de nombres de facteurs directs
}\endhead

L'objet de ce num\'ero est d'apporter dans certains cas
une pr\'ecision quantitative aux propri\'et\'es Min-fd et Max-fd
d\'efinies au num\'ero 6.

\medskip

   Soit $\Gamma$ un groupe. Pour un entier $n \ge 1$, 
notons $K_n(\Gamma)$ l'intersection de tous les sous-groupes normaux
de $\Gamma$ d'indices au plus $n$ 
et $k_n(\Gamma)$ l'indice de $K_n(\Gamma)$ dans $\Gamma$.

\medskip

   (i) Pour un produit direct, nous avons
$K_n(\prod_{j=1}^m \Gamma_j) = \prod_{j=1}^m K_n(\Gamma_j)$
et $k_n(\prod_{j=1}^m \Gamma_j) = \prod_{j=1}^m k_n(\Gamma_j)$.
Par exemple~: $k_n(\bold Z^m) = n^m$ pour tout $m \ge 1$.

Si $\Gamma$ est simple infini, alors $k_n(\Gamma) = 1$ pour tout $n \ge 1$.

\medskip

   (ii) Soit $\pi : \Gamma \longrightarrow \Delta$ un epimorphisme.
Alors $K_n(\Gamma) \subset \pi^{-1}(K_n(\Delta))$ et $k_n(\Gamma) \ge k_n(\Delta)$.
En particulier $k_n(\Gamma) \le k_n(F_g)$
pour un groupe $\Gamma$ \`a $g$ g\'en\'erateurs.

Pour une majoration grossi\`ere de $k_n(F_g)$, 
notons que $\log_n (k_n(F_g)) \le s_n^{\vartriangleleft}(F_g)$,
et rappelons que les r\'esultats du chapitre 2 de \cite{LuSe}
fournissent une majoration du nombre $s_n^{\vartriangleleft}(F_g)$ des sous-groupes
normaux d'indices au plus $n$ dans $F_g$.

\medskip

   (iii) Soit $\Gamma = \Gamma_1 \times \cdots \times \Gamma_m$~;
supposons que $\Gamma$ poss\`ede un syst\`eme de $g$ g\'en\'erateurs.
Notons $m_0$ le nombre des indices $j \in \{1 , \hdots, m\}$
tels que $k_n(\Gamma_j) \ge 2$. 
Il r\'esulte imm\'ediatement des points (i) et (ii) ci-dessus
que $m_0 \le \log_2(k_n(F_g))$.

\medskip

\proclaim{Th\'eor\`eme 13}
Soit $\Gamma$ un groupe qui poss\`ede un syst\`eme de $g$ g\'en\'erateurs
et qui est un produit $\Gamma = \Gamma_1 \times \cdots \times \Gamma_m$,
de groupes non r\'eduits \`a $\{1\}$.
S'il existe un entier $n \ge 2$ tel que $\Gamma$ est $n$-QM, 
alors $m \le \log_2(k_{n}(F_g))$.

En particulier, si $\Gamma$ est un groupe lin\'eaire de type fini, 
$\Gamma$ peut toujours s'\'ecrire comme produit direct 
d'un nombre fini de groupes ind\'ecomposables et
il existe une borne sur le nombre de facteurs 
des d\'ecompositions de $\Gamma$ en produits directs.
\endproclaim

\demo{D\'emonstration} Si le groupe de type fini $\Gamma$
est $n$-QM, chaque facteur $\Gamma_j$ l'est ausssi,
et la premi\`ere assertion de la proposition 
r\'esulte du point (iii) ci-dessus.
La seconde assertion r\'esulte alors de la proposition 12.
\hfill $\square$
\enddemo

\demo{Remarque} Etant donn\'e deux entiers $d,g \ge 1$ et un anneau de
type fini $A$, il r\'esulte des preuves ci-dessus qu'il existe une constante 
$M = M(d,g,A)$ telle que 
tout sous-groupe $\Gamma \subset \operatorname{GL}_d(A)$ \`a au plus $g$
g\'en\'erateurs poss\`ede une d\'ecomposition 
$\Gamma = \Gamma_1 \times \cdots \times \Gamma_m$ 
en produits de groupes ind\'ecomposables, avec $m \le M$.
\enddemo

Voici pour terminer une cons\'equence de la proposition  2 et du th\'eor\`eme 13.

\medskip

\proclaim{Proposition 14}
Soit $\Gamma$ un groupe lin\'eaire de type fini dont le centre est r\'eduit \`a un
\'el\'ement. Alors $\Gamma$ est uniquement directement d\'ecomposable.
\endproclaim

\medskip

Nous remercions Yves Benoist, Martin Bridson, Ken Brown, Luis Paris,  Alain Valette
et Thierry Vust  pour plusieurs observations pr\'ecieuses concernant notre texte.

\enddocument